\documentclass[a4paper,12pt,reqno]{amsart}  

\usepackage{amssymb,amsmath,amsfonts,amsthm,mathtools}
\usepackage{latexsym,mathrsfs,pb-diagram}
\usepackage[pdftex]{graphicx}
\usepackage{pgf,tikz,tikz-cd}

\usepackage[hmargin=3.0cm,vmargin=3.0cm]{geometry} 

\usepackage[plainpages=false,colorlinks,pdfpagelabels]{hyperref}
\hypersetup{
  urlcolor=black,
  citecolor=green,
  linkcolor=blue
}

\numberwithin{equation}{section}

\theoremstyle{definition}
\newtheorem{Def}{Definition}[section]

\theoremstyle{remark}

\newtheorem{Rem}[Def]{Remark}

\theoremstyle{plain}
\newtheorem{Prop}[Def]{Proposition}
\newtheorem{Cor}[Def]{Corollary}
\newtheorem{Thm}[Def]{Theorem}
\newtheorem{Lem}[Def]{Lemma}

\newcommand{\dfn}{\doteq}
\newcommand{\st}{ \ ; \ }
\newcommand{\lra}{\longrightarrow}
\newcommand{\sset}{\subset}
\newcommand{\bb}{\bullet}

\newcommand{\eset}{\emptyset}
\newcommand{\Z}{\mathbb{Z}}
\newcommand{\N}{\mathbb{N}}
\newcommand{\R}{\mathbb{R}}
\newcommand{\Q}{\mathbb{Q}}
\newcommand{\C}{\mathbb{C}}

\newcommand{\transp}[1]{\prescript{\mathrm{t} \!}{}{{#1}}}
\newcommand{\ort}{\mathrm{o}}

\DeclareMathOperator{\ran}{\mathrm{ran}}

\newcommand{\del}{\partial}
\newcommand{\dd}{\mathrm{d}}
\renewcommand{\Re}{\mathsf{Re}}
\renewcommand{\Im}{\mathsf{Im}}
\newcommand{\TT}{\mathbb{T}}

\newcommand{\D}{\mathscr{D}}

\newcommand{\cinfty}{\mathscr{C}^\infty}

\newcommand{\sob}{\mathscr{H}}
\newcommand{\sol}{\mathscr{S}}

\newcommand{\HH}{H}

\newcommand{\DD}{\mathrm{D}}
\newcommand{\VV}{\mathcal{V}}
\newcommand{\psum}{\sideset{}{'}\sum}

\newcommand{\LLambda}{\underline{\Lambda}}

\newcommand{\vv}{\mathrm}

\author{Gabriel Ara\'{u}jo}
\address{Universidade de S{\~a}o Paulo, Brazil}
\email{\texttt{gccsa@icmc.usp.br}}

\author{Igor A.~Ferra}
\address{Universidade Federal do ABC, Brazil}
\email{\texttt{ferra.igor@ufabc.edu.br}}

\author{Max R.~Jahnke}
\address{Philipps-Universit{\"a}t Marburg, Germany}
\email{\texttt{jahnkem@mathematik.uni-marburg.de}}

\author{Luis F.~Ragognette}
\address{Universidade Federal de Minas Gerais, Brazil}
\email{\texttt{luisragognette@mat.ufmg.br}}

\thanks{Supported by
  the S{\~a}o Paulo Research Foundation
  (FAPESP, grants~2016/13620-5, 2018/12273-5 and~2019/09967-8),
  by
  Conselho Nacional de Desenvolvimento Cient{\'i}fico e Tecnol{\'o}gico
  (CNPq, grant~163837/2022-8)
  and by
  Deutsche Forschungsgemeinschaft
  (DFG, grant~JA~3453/1-1).
}


\title[]{Global solvability and cohomology of tube structures on compact manifolds}

\begin{document}

\begin{abstract}
  We introduce new techniques
  to study the differential complexes
  associated to tube structures on
  $M \times \TT^m$ of corank $m$,
  in which
  $M$ is a compact manifold and
  $\TT^m$ is the $m$-torus.
  By systematically employing
  partial Fourier series,
  for complex tube structures,
  we completely characterize
  global solvability, in a given degree,
  in terms of a weak form of hypoellipticity,
  thus generalizing existing results
  and providing a broad answer to
  an open problem proposed by
  Hounie and Zugliani (2017).
  We also obtain new results
  on the finiteness of the
  cohomology spaces
  in intermediate degrees.
  In the case of real tube structures,
  we extend an isomorphism
  for the cohomology spaces
  originally obtained by
  Dattori da Silva and Meziani (2016)
  in the case
  $M = \TT^n$.
  Moreover, we establish
  necessary and sufficient conditions
  for the differential operator
  to have closed range
  in the first degree.
\end{abstract}

\maketitle


\section{Introduction}
\label{sec:defns}

We investigate the
global solvability and
the cohomology spaces of 
differential complexes associated with
certain systems of first-order differential operators
on a compact manifold,
the so-called
\emph{tube structures}.
These structures are invariant
under the action of a torus,
which allows us to use
partial Fourier series
to slice the non-elliptic differential complex
into infinitely many elliptic ones.
We can then apply
methods from elliptic theory
and sheaf cohomology
to study the resulting complexes.

Let
$M$ be a smooth, compact,
connected and oriented
$n$-dimensional manifold,
and let
$\omega_1, \ldots, \omega_m$ be smooth, complex, closed
$1$-forms on $M$,
and let
$\TT^m \dfn \R^m / 2 \pi \Z^m$.
On
$\Omega \dfn M \times \TT^m$,
we consider the involutive subbundle
$\VV \sset \C T \Omega$
whose sections are annihilated by 
\begin{equation}
  \label{eq:zetas}
  \zeta_k \dfn \dd x_k - \omega_k,
  \quad k \in \{1, \ldots, m\},
\end{equation}
in which
$x = (x_1, \ldots, x_m)$
denote the usual coordinates on
$\TT^m$.
Such
$\VV$ gives rise to a
complex of vector bundles and
first-order differential operators over
$\Omega$~\cite{treves_has, bch_iis},
here constructed as follows:
for each
$q \in \{0, \ldots, n\}$,
let
$\Lambda^q$ denote the bundle of
$q$-forms over $M$ and
$\LLambda^q$ its pullback via the projection
$\Omega \to M$. The
smooth sections of the latter
are locally written as
\begin{equation}
  \label{eq:form-coordinates}
  f = \psum_{|J| = q} f_J(t,x) \dd t_J,
\end{equation}
in which
$f_J \in \cinfty(U \times \TT^m)$
and
$(U; t_1, \ldots, t_n)$ is some local chart of $M$.
We define a differential operator
\begin{equation}
  \label{eq:def_dprime}
  \dd'
  \dfn
  \dd_t + \sum_{k = 1}^m \omega_k \wedge \frac{\del}{\del x_k}
  :
  \cinfty(\Omega; \LLambda^{q})
  \lra
  \cinfty(\Omega; \LLambda^{q + 1}),
\end{equation}
that satisfies
$\dd' \circ \dd' = 0$.
Our goal is to investigate
global solvability
-- that is,
\emph{closedness of the range} --
of~\eqref{eq:def_dprime}
in any degree
$q \in \{0, \ldots, n - 1\}$,
and to provide a better understanding of
the smooth cohomology spaces
\begin{equation}
  \label{eq:cohomology_VV}
  \HH^{q}_{\dd'}(\cinfty(\Omega))
  \dfn
  \frac{ \ker \{ \dd': \cinfty(\Omega; \LLambda^{q}) \lra \cinfty(\Omega; \LLambda^{q + 1}) \} }
       { \ran \{ \dd': \cinfty(\Omega; \LLambda^{q - 1}) \lra \cinfty(\Omega; \LLambda^{q}) \} },
\end{equation}
for
$q \in \{0, \ldots, n \}$
(with standard conventions
regarding the endpoints\footnote{All differential complexes
presented in this paper are
``completed as zero'' in negative degrees.}).

Our first result
is an equivalence between
global solvability of $\dd'$
and a weak notion of hypoellipticity 
(Theorem~\ref{thm:main_agh}),
a generalization
of~\cite[Corollary~7.2]{hz19}
that holds in arbitrary degree
and requires no further hypotheses
(except for being a tube),
and thus provides a broad answer
to the Open Problem 2
in~\cite{hz17}.
See also
Remark~\ref{rem:compatibility_conditions},
Corollary~\ref{cor:main_agh}
and
Remark~\ref{rem:hz_case}.

For real structures,
we prove certain
isomorphisms~\eqref{IsomorphismSobreGammaomega}
for the cohomology spaces,
similar to the main result
in~\cite{dm16}.
In that work,
$M$ is the $n$-torus,
and their result
roughly states that,
under suitable conditions, we have
\begin{equation}
  \label{eq:dm_iso_intro}
  \HH^{q}_{\dd'}(\cinfty(\Omega))
  \cong
  \cinfty(\TT^r) \otimes H^q_{\mathrm{dR}}(M),
\end{equation}
in which
$r$ is the rank of a group
associated to
$\omega_1, \ldots, \omega_m$.
For a general
$M$, however,
our approach reveals a new phenomenon:
an unexpected
obstruction~\eqref{eq:vanishing_mysterious_cohomologies}
for the validity
of~\eqref{eq:dm_iso_intro}
(Theorem~\ref{thm:strong_isomorphism}).
For surfaces,
this obstruction
-- which is likely related
both to the nature of the tube structure
and the topology of $M$ --
is not present
in any degree when
$M$ is either a
$2$-torus or a $2$-sphere,
but is present at
$q = 1$ when
$M$ has genus
$\mathsf{g} \geq 2$.
We show in
Section~\ref{exa:dm_fails}
how to construct tube structures
that satisfy all the hypothesis required
by~\cite{dm16},
but that do not
satisfy~\eqref{eq:dm_iso_intro}.
We conjecture
that~\eqref{eq:vanishing_mysterious_cohomologies}
holds in any degree when
$M$ is a compact Lie group
and the tube structure is real,
thus providing a complete generalization
of the results
in~\cite{dm16}.
We plan to
investigate this topic further
in a forthcoming work.

Still in the context of real structures,
we characterize the global solvability
of~\eqref{eq:def_dprime}
in the first degree
(Theorem~\ref{Thm:solv-diop})
via Diophantine conditions,
generalizing some results in~\cite{bcm93, bcp96}
for arbitrary corank
(thus a smooth version
of some results
in~\cite{adv23}).
Besides being more general,
our methods have advantages.
For example, the systematic use of the key
Lemma~\ref{lem:formal_closure} and its implications
results into relatively short
and straightforward proofs
without the need of usual techniques,
such as dualization with the top degree,
or the use of
\emph{a priori} inequalities. 

In Sections~\ref{sec:fourier_cohom_spaces}
and~\ref{sec:gen_fin_thm},
we revisit complex tube structures,
and derive a handful of necessary conditions
for finiteness of the cohomology spaces. 
The results we proved
led us to conjecture whether
$\HH^{q}_{\dd'}(\cinfty(\Omega))$
is finite dimensional
only if
it is isomorphic with
$H^q_{\mathrm{dR}}(M)$
-- true, notably,
when~\eqref{eq:vanishing_mysterious_cohomologies}
holds and $r = 0$ 
(Theorem~\ref{thm:finiteness_when_mysterious}),
as well as in every other situation
in which we were able to compute.
For example, when
$M$ is a surface, we obtain
a quite complete description
(Section~\ref{sec:surfaces}).

\section{Preliminaries}

\subsection{Global solvability
  in abstract complexes
  and related concepts}

Let
$X$ be a smooth, compact,
connected and oriented manifold,
and
$\mathbb{E}$
be a complex vector bundle over $X$.
The space
$\cinfty(X; \mathbb{E})$
of smooth sections of
$\mathbb{E}$
carries its standard
Fr\'echet topology;
by endowing
$X$ with a Riemannian metric and
$\mathbb{E}$ with a Hermitian metric,
one may write it
as the projective limit
of a suitable sequence
of Sobolev spaces of sections of
$\mathbb{E}$.
We also let
$\D'(X; \mathbb{E})$
be the space of distribution sections of
$\mathbb{E}$,
which will be identified
with the topological dual of
$\cinfty(X; \mathbb{E}^*)$.

Let
$\mathbb{E}, \mathbb{F}$
be vector bundles over
$X$ and
$P$ a differential operator from
$\mathbb{E}$ to $\mathbb{F}$.

\begin{Def}
  We say that $P$ is
  \emph{almost globally hypoelliptic}
  $\mathrm{(AGH)}$
  if
  \begin{equation*}
    \forall u \in \D'(X; \mathbb{E}), \ Pu \in \cinfty(X; \mathbb{F})
    \Longrightarrow
    \text{
      $\exists v \in \cinfty(X; \mathbb{E})$
      such that
      $Pv = Pu$
    }. 
  \end{equation*}
\end{Def}
We have:
\begin{Thm}
  \label{thm:agh_gs}
  If
  $P$ is~$\mathrm{(AGH)}$,
  then
  $P: \cinfty(X; \mathbb{E}) \to \cinfty(X; \mathbb{F})$
  has closed range.
\end{Thm}
This result is proved 
in~\cite[Theorem~3.5]{afr22}
for scalar operators;
its proof extends
to vector-valued operators
in a straightforward way,
hence we omit it.
A converse fails to hold
even for very simple operators,
but is valid
for many classes of operators~\cite{araujo19, afr22},
including~\eqref{eq:def_dprime},
as we will prove in
Theorem~\ref{thm:main_agh}
below.

Let
$\mathbb{G}$ be a third vector bundle over
$X$ and
$Q$ a differential operator from
$\mathbb{F}$ to $\mathbb{G}$
such that
$Q \circ P = 0$.
We define two cohomology spaces
\begin{equation*}
  \HH_{P,Q}(\mathscr{F}(X))
  \dfn
  \frac{\ker \{ Q: \mathscr{F}(X; \mathbb{F}) \lra \mathscr{F}(X; \mathbb{G}) \}}
       {\ran \{ P: \mathscr{F}(X; \mathbb{E}) \lra \mathscr{F}(X; \mathbb{F}) \}},
  \quad \mathscr{F} = \text{$\cinfty$ or $\D'$}.
\end{equation*}
Our goal is to understand
these two cohomology spaces separately,
as well as their relationship.
The inclusions
$\cinfty(X; \ast) \hookrightarrow \D'(X; \ast)$
($\ast = \mathbb{E}, \mathbb{F}, \mathbb{G}$)
induce a linear morphism
\begin{equation}
  \label{eq:cinftytodistrcohomology}
  \HH_{P,Q}(\cinfty(X)) \lra \HH_{P,Q}(\D'(X)) 
\end{equation}
which may be neither
injective nor surjective.
However,
\begin{equation*}
  \text{$P$ is~$\mathrm{(AGH)}$}
  \Longleftrightarrow
  \text{\eqref{eq:cinftytodistrcohomology}
    is injective},
\end{equation*}
being therefore
independent of $Q$,
and a property of $P$ alone.
One should recall the
traditional notion of
global hypoellipticity
for complexes,
that is:
\begin{multline}
  \label{eq:GHforsystems}
  \forall f \in \D'(X; \mathbb{F}), \ Qf \in \cinfty(X; \mathbb{G})
  \Longrightarrow \\
  \text{
    $\exists u \in \D'(X; \mathbb{E})$
    such that
    $f - Pu \in \cinfty(X; \mathbb{F}).$
  }
\end{multline}
Since
$Q (f - Pu) = Qf$,~\eqref{eq:GHforsystems} implies that
$Q$ is~$\mathrm{(AGH)}$. Additionally, we have that
\begin{equation*}
  \forall f \in \D'(X; \mathbb{F}), \ Qf = 0
  \Longrightarrow
  \text{
    $\exists u \in \D'(X; \mathbb{E})$
    such that
    $f - Pu \in \cinfty(X; \mathbb{F}),$
  }
\end{equation*}
which is equivalent
to the surjectivity
of~\eqref{eq:cinftytodistrcohomology}.
We state this result more precisely:
\begin{Prop}
  Property~\eqref{eq:GHforsystems} holds
  if and only if
  $Q$ is~$\mathrm{(AGH)}$ 
  and~\eqref{eq:cinftytodistrcohomology}
  is onto.
\end{Prop}

The transpose
$\transp{P}$ 
is a differential operator from
$\mathbb{F}^*$ to $\mathbb{E}^*$,
yielding new maps
\begin{equation*}
  \transp{P}: \left\{
  \begin{array}{c c c}
    \D'(X; \mathbb{F}^*)     & \lra & \D'(X; \mathbb{E}^*) \\
    \cinfty(X; \mathbb{F}^*) & \lra & \cinfty(X; \mathbb{E}^*)
  \end{array},
  \right.
\end{equation*}
the latter being
the restriction of the former.
In the presence of
a second operator $Q$
satisfying $Q \circ P = 0$,
we want to determine
necessary conditions on
$f \in \cinfty(X; \mathbb{F})$
so as to be able to solve
$Pu = f$ with
$u \in \cinfty(X; \mathbb{E})$.
Obviously, we must have
\begin{equation}
  \label{eq:local_cc}
  f \in \ker \{ Q: \cinfty(X; \mathbb{F}) \lra \cinfty(X; \mathbb{G}) \}. 
\end{equation}
Moreover, if
$v \in \D'(X; \mathbb{F}^*)$
is such that
$\transp{P}v = 0$,
then
\begin{equation*}
  \langle v, f \rangle
  = \langle v, Pu \rangle = \langle \transp{P}v, u \rangle
  = 0,
\end{equation*}
that is,
\begin{equation}
  \label{eq:global_cc}
  f \in \ker \{ \transp{P}: \D'(X; \mathbb{F}^*) \lra \D'(X; \mathbb{E}^*) \}^\ort. 
\end{equation}
However,
by Functional Analysis,
the annihilator in~\eqref{eq:global_cc}
equals the closure of
$\ran \{ P: \cinfty(X; \mathbb{E}) \to \cinfty(X; \mathbb{F}) \}$,
and since
$Q \circ P = 0$,
the range of $P$ is contained in
$\ker \{ Q: \cinfty(X; \mathbb{F}) \to \cinfty(X; \mathbb{G}) \}$,
which is closed in
$\cinfty(X; \mathbb{F})$.
We conclude:
\begin{equation*}
  \ker \{ \transp{P}: \D'(X; \mathbb{F}^*) \lra \D'(X; \mathbb{E}^*) \}^\ort
  \sset
  \ker \{ Q: \cinfty(X; \mathbb{F}) \lra \cinfty(X; \mathbb{G}) \},
\end{equation*}
hence 
the compatibility condition~\eqref{eq:local_cc}
is redundant in light
of~\eqref{eq:global_cc}.
We therefore introduce
the following definition
in spite of the presence of the operator $Q$.
\begin{Def}
  \label{def:solvability_compact}
  We say that $P$ is
  \emph{globally solvable}
  if for every
  $f \in \cinfty(X; \mathbb{F})$
  satisfying
  \begin{equation*}
    \langle v, f \rangle = 0
    \quad
    \text{for every
      $v \in \D'(X; \mathbb{F}^*)$
      such that
      $\transp{P}v = 0$
    } 
  \end{equation*}
  there exists
  $u \in \cinfty(X; \mathbb{E})$
  such that $Pu = f$.
\end{Def}
In other words,
$P$ is globally solvable
if and only if
the range of
$P$ is closed.

\subsection{Partial Fourier series
  for  sections of $\LLambda^{q}$}
\label{sec:partial_fourier}

Given
$f \in \cinfty(\Omega; \LLambda^{q})$,
for each $\xi \in \Z^m$, we define
an element
$\hat{f}_\xi \in \cinfty(M; \Lambda^q)$
as follows:
if $U \sset M$
is a coordinate open set
in which $f$ is written
as~\eqref{eq:form-coordinates}
then
\begin{equation*}
  \hat{f}_\xi(t)
  \dfn \psum_{|J| = q} \hat{f}_J(t,\xi) \dd t_J
  = \psum_{|J| = q}
  \left( \int_{\TT^m} e^{-i x \xi} f_J(t, x) \dd x \right) \dd t_J,
  \quad t \in U.
\end{equation*}
Note that since
$\LLambda_{(t,x)}= \Lambda_t$
for every
$x \in \TT^m$,
each
$\dd t_J$
can be thought
as a section of either
$\Lambda^q$ or
$\LLambda^q$.
One can check that
the construction above
is independent of
the choice of coordinates on $U$,
hence defines well
a differential form
of degree $q$ in $M$.

Indeed, suppose 
$\chi_1, \ldots, \chi_n$
is a frame of complex $1$-forms
across $U$;
we denote by
$\vv{X}_1, \ldots, \vv{X}_n \in \cinfty(U; \C T M)$
the associated dual frame.
Then
\begin{equation*}
  \chi_J \dfn \chi_{j_1} \wedge \cdots \wedge \chi_{j_q},
  \quad \text{for $J = (j_1, \ldots, j_q) \in \{1, \ldots, n\}^q$ ordered},
\end{equation*}
form a frame for $\Lambda^q$ on $U$.
We may write
$f$
on
$U \times \TT^m$ as
\begin{equation}
  \label{eq:form_frames}
  f = \psum_{|J| = q} f_J(t,x) \chi_J,
  \quad f_J \dfn f(\vv{X}_J) \in \cinfty(U \times \TT^m).
\end{equation}
If
$\chi_1^\bb, \ldots, \chi_n^\bb$
is another coframe in $U$,
with dual frame
$\vv{X}_1^\bb, \ldots, \vv{X}_n^\bb$,
then
\begin{equation*}
  \vv{X}_j = \sum_{k = 1}^n a_{j k} \vv{X}_{k}^\bb,
  \quad a_{jk} \dfn \langle \chi_k^\bb, \vv{X}_j \rangle \in \cinfty(U),
\end{equation*}
and more generally
for higher multi-indices
\begin{equation*}
  \vv{X}_J = \psum_{|K| = q} a_{JK} \vv{X}_{K}^\bb,
  \quad a_{JK} \dfn \langle \chi_K^\bb, \vv{X}_J \rangle \in \cinfty(U).
\end{equation*}
Hence in~\eqref{eq:form_frames}
we have
\begin{equation*}
  f_J = f(\vv{X}_J)
  = \psum_{|K| = q} a_{JK} f(\vv{X}_K^\bb)
  = \psum_{|K| = q} a_{JK} f_K^\bb
\end{equation*}
and noticing that
$a_{JK}$ does not depend on $x$
we have
\begin{equation*}
  \hat{f}_J(\cdot,\xi)
  = \int_{\TT^m} e^{-i x \xi} f_J(\cdot, x) \dd x
  = \psum_{|K| = q} a_{JK} \int_{\TT^m} e^{-i x \xi} f_K^\bb (\cdot, x) \dd x
  = \psum_{|K| = q} a_{JK} \hat{f}_K^\bb (\cdot,\xi)
\end{equation*}
therefore
\begin{equation*}
  \psum_{|J| = q} \hat{f}_J(\cdot,\xi) \chi_J
  = \psum_{|K| = q} \hat{f}_K^\bb (\cdot,\xi) \psum_{|J| = q} a_{J K}  \chi_J
  = \psum_{|K| = q} \hat{f}_K^\bb (\cdot,\xi) \chi_K^\bb.
\end{equation*}

It is useful
to regard
$\hat{f}_\xi$
as a current on $M$, that is, we define
$\hat{f}_\xi: \cinfty(M; \Lambda^{n - q}) \to \C$ by
\begin{equation}
  \label{eq:fourier_by_duality}
 g \mapsto \langle \hat{f}_\xi, g \rangle_M
  \dfn \langle f, e^{-ix \xi} \wedge g \wedge \dd x \rangle
  = \int_\Omega f \wedge e^{-ix \xi} \wedge g \wedge \dd x, 
\end{equation}
with
$\dd x \dfn \dd x_1 \wedge \cdots \wedge \dd x_m$;
both definitions yield the same object.
Indeed,
we fix a Riemannian metric on
$M$ and
evaluate~\eqref{eq:fourier_by_duality}
on a
$g \in \cinfty_c(U; \Lambda^{n - q})$,
in which
$U \sset M$
is the domain of an orthonormal coframe
$\chi_1, \ldots, \chi_n$,
which we write
\begin{equation*}
  g = \psum_{|K| = n - q} g_K(t) \chi_K,
  \quad g_K \in \cinfty_c(U).
\end{equation*}
Each ordered 
$J \in \{1, \ldots, n\}^q$
pairs with a single ordered
$J^\star \in \{1, \ldots, n\}^{n - q}$
such that
\begin{equation}
  \label{eq:Jstar}
  \chi_J \wedge \chi_{J^\star}
  = \epsilon_J \ \chi_1 \wedge \cdots \wedge \chi_n
  = \epsilon_J \ \dd V_M,
  \quad \epsilon_J = \pm 1,
\end{equation}
and
$\chi_J \wedge \chi_K = 0$
for the remaining $K$
(we denote by $\dd V_M$
the Riemannian volume form on $M$).
Therefore,
using~\eqref{eq:form_frames},
\begin{equation*}
  \int_\Omega f \wedge e^{-ix \xi} \wedge g \wedge \dd x
  = \psum_{|J| = q} \int_{U \times \TT^m} e^{-i x \xi} f_J(t,x) g_{J^\star}(t) \epsilon_J \dd V_M(t) \wedge \dd x
\end{equation*}
while according to
the earlier local definition
\begin{equation*}
  \int_M \hat{f}_\xi \wedge g
  = \psum_{|J| = q} \int_U \left( \int_{\TT^m} e^{-i x \xi} f_J(t, x) \dd x \right) g_{J^\star}(t) \epsilon_J \dd V_M(t)
\end{equation*}
which certainly yields
the same result.
This extends our construction to
$q$-currents,
yielding linear maps
\begin{equation*}
  \mathcal{F}_\xi :
  \left\{
  \begin{array}{c c c}
    \cinfty(\Omega; \LLambda^{q}) & \lra & \cinfty(M; \Lambda^q) \\
    \D'(\Omega; \LLambda^{q})     & \lra & \D'(M; \Lambda^q)
  \end{array}
  \right.
\end{equation*}
defined by the assignment
$f \mapsto \hat{f}_\xi$.

\begin{Prop}
  \label{prop:continuity_Fxi}
  The map
  $\mathcal{F}_\xi: \cinfty(\Omega; \LLambda^{q}) \to \cinfty(M; \Lambda^q)$
  is continuous.
\end{Prop}
\begin{proof}
  Since
  $\mathcal{F}_\xi$ is linear,
  it is enough to prove that given 
  $\{ f_\nu \}_{\nu \in \N} \sset \cinfty(\Omega; \LLambda^{q})$
  converging to zero 
  we have 
  $\mathcal{F}_\xi(f_\nu) \to 0$
  in
  $\cinfty(M; \Lambda^q)$.
  Taking
  $U_1, \ldots, U_\ell$ a finite
  coordinate open cover of
  $M$,
  we can write
  \begin{equation*}
     f_\nu|_{U_p\times \TT^m}   = \psum_{|J| = q} f_{\nu,p ,J}(t,x) \dd t_J,
  \end{equation*}
  in which
  $f_{\nu,p, J} \in \cinfty(U_p \times \TT^m)$
  for every
  $p \in \{1, \ldots, \ell\}$ and
  $|J| = q$.
  We can assume that each
  $f_{\nu,p, J}$ is the restriction of
  a smooth function defined
  on an open subset of
  $M \times \TT^m$ that contains
  $\overline{U_p} \times \TT^m$,
  hence our hypothesis on
  $\{ f_\nu \}_{\nu \in \N}$
  implies that given
  $\epsilon > 0$ we have,
  for every $\alpha \in \Z_+^n$,
  that there exists
  $\nu_\alpha\in \N$
  such that
  for every
  $\nu \geq \nu_\alpha$,
  $p \in \{ 1, \ldots, \ell \}$
  and
  $|J| = q$
  we have
  \begin{equation*}
    \sup_{(t,x) \in U_p \times \TT^m}  \| \del_t^{\alpha} f_{\nu, p, J}(t,x) \| < \epsilon
  \end{equation*}
  and, consequently,
  \begin{equation*}
    \sup_{t\in U_p} \| \del^{\alpha} \mathcal{F}_\xi(f_\nu)(t) \|
    =
    \sup_{t\in U_p} \Big\| \psum_{|J| = q} \left( \int_{\TT^m} e^{-i x \xi} \del^{\alpha}_tf_{\nu,p, J}(t, x) \dd x \right) \dd t_J \Big\|
    \leq
    (2\pi)^m C \epsilon
  \end{equation*}
  for some
  $C > 0$.
  This shows that
  $\mathcal{F}_\xi(f_\nu) \to 0$ in
  $\cinfty(M; \Lambda^q)$,
  thus proving the continuity of
  $\mathcal{F}_\xi$.
 \end{proof}

Going in the other direction
we define linear maps
\begin{equation}
  \label{eq:Exi}
  \mathcal{E}_\xi :
  \left\{
  \begin{array}{c c c}
    \cinfty(M; \Lambda^q) & \lra & \cinfty(\Omega; \LLambda^{q}) \\
    \D'(M; \Lambda^q)     & \lra & \D'(\Omega; \LLambda^{q})
  \end{array}
  \right.
\end{equation}
by
$f \mapsto (2 \pi)^{-m} e^{ix \xi} \wedge f$.
\begin{Lem}
  \label{lem:rightinverseFxiExi}
  For every
  $\xi \in \Z^m$,
  we have that
  $\mathcal{F}_\xi \circ \mathcal{E}_\xi$
  is the identity on
  $\cinfty(M; \Lambda^q)$
  and
  $\mathcal{F}_\xi \circ \mathcal{E}_\eta = 0$ if
  $\eta \neq \xi$,
  for every
  $q \in \{0,\ldots, n\}$.
\end{Lem}
\begin{proof}
  Let
  $f \in \cinfty(M; \Lambda^q)$.
  Given
  $g \in \cinfty(M; \Lambda^{n - q})$
  and
  $\xi, \eta \in \Z^m$
  we have
  \begin{equation*}
    \langle \mathcal{F}_\xi \mathcal{E}_\eta f, g \rangle_M
    = \frac{1}{(2 \pi)^m} \int_\Omega e^{ix\eta} \wedge f \wedge e^{-ix\xi} \wedge g \wedge \dd x
    = \frac{1}{(2 \pi)^m} \int_{\TT^m} e^{ix(\eta - \xi)} \dd x \int_M f \wedge g
  \end{equation*}
  hence
  $\mathcal{F}_\xi \mathcal{E}_\eta f = \delta_{\xi \eta} f$.
\end{proof}

\begin{Prop}
  Given
  $f \in \cinfty(\Omega; \LLambda^{q})$
  we have 
  \begin{equation}
    \label{eq:THE_fourier_series}
    f = \frac{1}{(2 \pi)^m} \sum_{\xi \in \Z^m} e^{ix \xi} \hat{f}_\xi
  \end{equation}
  with convergence in
  $\cinfty(\Omega; \LLambda^{q})$.
  In particular,
  $f = 0$
  if and only if
  $\hat{f}_\xi = 0$
  for all
  $\xi \in \Z^m$.
\end{Prop}
\begin{proof}
  Note that if
  $\chi \in \cinfty(M)$
  then
  $\widehat{(\chi f)}_\xi = \chi \hat{f}_\xi$.
  Thus, it is enough to verify the
  convergence assuming that
  the support of
  $f$ is contained in
  $W \times \TT^m$, in which
  $W$ is a coordinate open subset of
  $M$.
  We can assume that
  there exists a diffeomorphism
  $\lambda: (-\pi, \pi)^{n} \to W$.
  Note that
  $\mathcal{F}_\xi(\lambda_\TT^* f) = \lambda^* \hat{f}_\xi$,
  in which
  $\lambda_\TT \dfn \lambda \times \mathrm{id}_{\TT^m}$.
  Since we can extend
  $\lambda_{\TT}^* f$ periodically to
  $\R^n\times \TT^{m}$
  and identify it with a smooth
  $q$-form
  $F$ defined
  on $\TT^n \times \TT^m$,
  we can use partial Fourier expansion
  on $\TT^n \times \TT^m$:
  \begin{multline*}
    \lambda_{\TT}^* f(t,x)
    =
    F(t,x)
    =
    \frac{1}{(2 \pi)^m} \sum_{\xi \in \Z^m} e^{ix\xi} \hat{F}_\xi(t)
    =
    \frac{1}{(2 \pi)^m} \sum_{\xi \in \Z^m} e^{ix\xi} \mathcal{F}_\xi (\lambda_{\TT}^* f)(t) \\
    =
    \frac{1}{(2 \pi)^m} \sum_{\xi \in \Z^m} e^{ix\xi}( \lambda^* \hat{f}_\xi)(t)
    =
    \lambda_{\TT}^* \bigg( \frac{1}{(2 \pi)^m} \sum_{\xi \in \Z^m} e^{ix\xi} \hat{f}_\xi \bigg)(t,x).
  \end{multline*}
  Since
  $\lambda^*_\TT$ is a
  topological isomorphism between
  the corresponding space of
  $q$-forms,
  we conclude the promised identity.
\end{proof}

\subsubsection{Hodge star and adjoints}

Assume again
$M$ endowed with a
Riemannian metric.
Given $t \in M$,
for each
$g_0 \in \Lambda^q_t$, there is a unique
$g^\star_0 \in \Lambda^{n - q}_t$
such that
\begin{equation*}
  f_0 \wedge \bar{g}^\star_0
  = \langle f_0, g_0 \rangle_{\Lambda^q_t} \ \dd V_M(t),
  \quad \forall f_0 \in \Lambda^q_t.
\end{equation*}
In terms of an orthonormal coframe
$\chi_1, \ldots, \chi_n$,
we write
\begin{equation*}
  f_0 = \psum_{|J| = q} f_{0,J} \chi_J,
  \quad
  g_0 = \psum_{|J| = q} g_{0,J} \chi_{J},
  \quad
  g_0^\star = \psum_{|J| = q} g_{0,J^\star}^\star \chi_{J^\star},
\end{equation*}
in which
$J^\star$ is defined in~\eqref{eq:Jstar}.
Hence
\begin{equation*}
  \langle f_0, g_0 \rangle_{\Lambda^q_t}
  = \psum_{|J| = q} f_{0,J} \bar{g}_{0,J},
  \quad
  f_0 \wedge \bar{g}^\star_0
  = \psum_{|J| = q} f_{0,J} \bar{g}_{0,J^\star}^\star \epsilon_J \ \dd V_M(t),
\end{equation*}
which allows us to find
$g_{0,J^\star}^\star \dfn \epsilon_J g_{0,J}$
for every $J$.
This gives rise to
a linear map
$\Lambda^q_t \to \Lambda^{n - q}_t$,
and hence to a bundle map
$\Lambda^q \to \Lambda^{n - q}$.
The induced map on sections
\begin{equation*}
  g \in \cinfty(M; \Lambda^q)
  \longmapsto
  g^\star \in \cinfty(M; \Lambda^{n - q})
\end{equation*}
satisfies
\begin{equation*}
  \langle f, g \rangle_{L^2(M; \Lambda^q)}
  = \int_M \langle f(t), g(t) \rangle_{\Lambda^q_t} \ \dd V_M(t)
  = \int_M f \wedge \bar{g}^\star,
  \quad \forall f \in \cinfty(M; \Lambda^q).
\end{equation*}
On $\Omega$,
we have a map on sections
\begin{equation*}
  g \in \cinfty(\Omega; \LLambda^{q})
  \longmapsto
  g^\star \in \cinfty(\Omega; \LLambda^{n - q})
\end{equation*}
such that for every
$f \in \cinfty(\Omega; \LLambda^{q})$:
\begin{align*}
  \langle f, g \rangle_{L^2(\Omega; \LLambda^{q})}
  &= \int_\Omega \langle f(t,x), g(t,x) \rangle_{\Lambda^q_t} \ \dd V(t,x) \\
  &= \int_{\TT^m} \int_M \langle f(t,x), g(t,x) \rangle_{\Lambda^q_t} \ \dd V_M(t) \dd x \\
  &= \int_{\TT^m} \langle f(\cdot,x), g(\cdot,x) \rangle_{L^2(M; \Lambda^q)} \dd x \\
  &= \int_{\TT^m} \int_M f(\cdot,x) \wedge \bar{g}^\star(\cdot,x) \dd x \\
  &= \int_\Omega f \wedge \bar{g}^\star \wedge \dd x,
\end{align*}
in which
$\dd V$ stands for
the volume form
associated with the product metric
on $\Omega$.

\begin{Lem}
  \label{lem:adjoints}
  The maps
  $\mathcal{F}_\xi$ and
  $(2 \pi)^m \mathcal{E}_{\xi}$
  are formal adjoints
  of each other.
\end{Lem}
\begin{proof}
  For
  $f \in \cinfty(\Omega; \LLambda^{q})$
  and
  $g \in \cinfty(M; \Lambda^q)$,
  we compute
  \begin{multline*}
    \langle \hat{f}_\xi, g \rangle_{L^2(M; \Lambda^q)}
    = \langle \hat{f}_\xi, \bar{g}^\star \rangle_M
    = \langle f, e^{-ix\xi}  \wedge \bar{g}^\star \wedge \dd x \rangle \\
    = \langle f, \overline{(e^{ix\xi} \wedge g)^\star} \wedge \dd x \rangle
    = \langle f, e^{ix\xi} \wedge g \rangle_{L^2(\Omega; \LLambda^{q})}.
  \end{multline*}
\end{proof}

\subsection{Fourier analysis of $\dd'$}
\label{sec:fourierL}

Let
$\pmb{\omega} \dfn (\omega_1, \ldots, \omega_m)$
with each
$\omega_k$ being a closed
$1$-form on $M$.
Given
$\xi \in \Z^m$,
we write
\begin{equation*}
  \xi \cdot \pmb{\omega} \dfn \sum_{k = 1}^m \xi_k \omega_k \in \cinfty(M; \Lambda^1)
\end{equation*}
and define
$\dd'_\xi: \cinfty(M; \Lambda^q) \to \cinfty(M; \Lambda^{q + 1})$
by
\begin{equation*}
  \dd'_\xi f \dfn \dd f + i (\xi \cdot \pmb{\omega}) \wedge f,
\end{equation*}
a first-order differential operator
that satisfies
$\dd'_\xi \circ \dd'_\xi = 0$, thus forming a complex whose
smooth cohomology spaces we denote by
\begin{equation}
  \label{eq:cohomology-spaces-dprimexi}
  \HH_\xi^q (\cinfty(M))
  \dfn
  \frac{\ker \{ \dd'_\xi: \cinfty(M; \Lambda^q) \lra \cinfty(M; \Lambda^{q + 1}) \}}
       {\ran \{ \dd'_\xi: \cinfty(M; \Lambda^{q - 1}) \lra \cinfty(M; \Lambda^q) \}},
  \quad q \in \{1, \ldots, n\}
\end{equation}
and, as usual,
$\HH_\xi^0(\cinfty(M)) \dfn \ker \{ \dd'_\xi: \cinfty(M) \lra \cinfty(M; \Lambda^{ 1}) \}$.

In the next section,
we start a deeper study
of these zero-order perturbations of the
de Rham complex.
For the moment,
we focus on their formal aspects
related to the Fourier analysis
of the complex $\dd'$.

\begin{Lem}
  \label{lem:transposes}
  The following transposition formulas hold.
  \begin{enumerate}
  \item For
    $f \in \cinfty(\Omega; \LLambda^{q})$
    and
    $g \in \cinfty(\Omega; \LLambda^{n - q - 1})$,
    we have
    \begin{equation*}
      \int_\Omega \dd' f \wedge g \wedge \dd x = (-1)^{q + 1} \int_\Omega f \wedge \dd' g \wedge \dd x.
    \end{equation*}
  \item For
    $f \in \cinfty(M; \Lambda^q)$
    and
    $g \in \cinfty(M; \Lambda^{n - q - 1})$,
    we have
    \begin{equation*}
      \int_M \dd'_\xi f \wedge g = (-1)^{q + 1} \int_M f \wedge \dd'_{-\xi} g
    \end{equation*}
    for every
    $\xi \in \Z^m$.
  \end{enumerate}
\end{Lem}
\begin{proof}
  \begin{enumerate}
  \item Notice that
    \begin{equation*}
      \dd f \wedge \zeta = \dd' f \wedge \zeta,
    \end{equation*}
    in which
    $\zeta \dfn \zeta_1 \wedge \cdots \wedge \zeta_m$;
    recall that
    $\zeta_k = \dd x_k - \omega_k$,
    hence
    $\zeta$ is a closed
    $m$-form on
    $\Omega$.
    It follows that
    \begin{align*}
      \dd (f \wedge g \wedge \zeta)
      &= \dd f \wedge g \wedge \zeta
      + (-1)^q f \wedge \dd g \wedge \zeta
      + (-1)^{n - 1} f \wedge g \wedge \dd \zeta \\
      &= \dd' f \wedge g \wedge \zeta
      + (-1)^q f \wedge \dd' g \wedge \zeta
    \end{align*}
    thanks to a twofold application
    of the previous argument,
    so by Stokes Theorem we have
    \begin{equation*}
      \int_\Omega \dd' f \wedge g \wedge \zeta
      = (-1)^{q + 1} \int_\Omega f \wedge \dd' g \wedge \zeta.
    \end{equation*}
    Since it is also clear that
    $h \wedge \zeta = h \wedge \dd x$
    whenever
    $h \in \cinfty(\Omega; \LLambda^{n})$,
    the conclusion follows at once.
  \item Similarly:
    \begin{align*}
      \int_M \dd'_\xi f \wedge g
      &= \int_M \dd f \wedge g + \int_M (i \xi \cdot \pmb{\omega}) \wedge f \wedge g \\
      &= (-1)^{q + 1} \int_M f \wedge \dd g + (-1)^q \int_M f \wedge (i \xi \cdot \pmb{\omega}) \wedge g \\
      &= (-1)^{q + 1} \int_M f \wedge \dd g + (-1)^{q + 1} \int_M f \wedge (-i \xi \cdot \pmb{\omega}) \wedge g \\
      &= (-1)^{q + 1} \int_M f \wedge \dd'_{-\xi} g.
    \end{align*}
  \end{enumerate}
\end{proof}

\begin{Cor}
  \label{cor:actions_fouriercoeff}
  For each
  $\xi \in \Z^m$,
  the following commutation relations hold:
  \begin{enumerate}
  \item $\dd' \circ \mathcal{E}_\xi = \mathcal{E}_\xi \circ \dd'_\xi$;
  \item $\mathcal{F}_\xi \circ \dd' = \dd'_\xi \circ \mathcal{F}_\xi$.
  \end{enumerate}
\end{Cor}
\begin{proof}
  \begin{enumerate}
  \item For
    $f \in \cinfty(M; \Lambda^q)$,
    we have
    \begin{align*}
      \dd' (e^{ix \xi} \wedge f)
      &= \dd_t (e^{ix \xi} \wedge f) + \sum_{k = 1}^m \omega_k \wedge \del_{x_k} (e^{ix \xi} \wedge f) \\
      &= e^{ix \xi} \wedge \dd_t f + \sum_{k = 1}^m \omega_k \wedge i\xi_k e^{ix \xi} \wedge f \\
      &= e^{ix \xi} \wedge \dd'_\xi f.
    \end{align*}
  \item Given
    $f \in \cinfty(\Omega; \LLambda^{q})$,
    for an arbitrary
    $g \in \cinfty(M; \Lambda^{n - q - 1})$,
    we can use
    the previous item and
    Lemmas~\ref{lem:adjoints}
    and~\ref{lem:transposes}
    to compute:
    \begin{align*}
      \langle \mathcal{F}_\xi( \dd' f), g \rangle_M
      &= (2 \pi)^m \langle \dd' f, (\mathcal{E}_{-\xi} g) \wedge \dd x \rangle \\
      &= (2 \pi)^m (-1)^{q + 1} \langle f, [\dd' (\mathcal{E}_{-\xi} g)] \wedge \dd x \rangle \\
      &= (2 \pi)^m (-1)^{q + 1} \langle f, [\mathcal{E}_{-\xi} (\dd'_{-\xi} g)] \wedge \dd x \rangle \\
      &= (-1)^{q + 1} \langle \mathcal{F}_\xi f, \dd'_{-\xi} g \rangle_M \\
      &= \langle \dd'_{\xi} (\mathcal{F}_\xi f), g \rangle_M.
    \end{align*}
    Hence
    $\mathcal{F}_\xi( \dd' f) = \dd'_{\xi} (\mathcal{F}_\xi f)$. \qedhere
  \end{enumerate}
\end{proof}

\section{Zero-order perturbations of the de Rham complex}
\label{sec:zero_order_pert1}

Motivated by
the introduction of the operators
$\dd'_\xi$ in
Section~\ref{sec:fourierL},
we discuss a general class
of operators on $M$
that are perturbations
of the exterior derivative.
Let
$\omega$ be a complex, smooth,
closed $1$-form on $M$,
and define,
for $q \in \{0, \ldots, n\}$:
\begin{equation}
  \label{eq:abstract_perturbation}
  \DD_\omega \dfn \dd + i \omega \wedge \cdot :
  \cinfty(M; \Lambda^q) \lra \cinfty(M; \Lambda^{q + 1}).
\end{equation}
Note that
$\dd'_\xi= \DD_{\omega}$ when
$\omega= \xi \cdot \pmb{\omega}$. 
These are first-order differential operators,
satisfying
$\DD_\omega \circ \DD_\omega = 0$
since
$\omega$ is closed.
The operators $\DD_\omega$ define a differential complex,
which is elliptic since
$\DD_\omega$
has the same principal symbol as $\dd$.
Such a differential complex,
however,
does not come from
an involutive structure on $M$
(for instance,
$\DD_\omega(1) = i\omega \neq 0$),
hence we cannot apply
solvability results from
this theory to it.

Given any open set
$U \sset M$, 
we define
\begin{equation}
  \label{eq:cohomology-spaces-Domega}
  \HH_\omega^q(\mathscr{F}(U)) \dfn
  \frac{\ker \{ \DD_\omega: \mathscr{F}(U; \Lambda^q) \lra \mathscr{F}(U; \Lambda^{q + 1}) \}}
       {\ran \{ \DD_\omega: \mathscr{F}(U; \Lambda^{q - 1}) \lra \mathscr{F}(U; \Lambda^q) \}},
  \quad \mathscr{F} = \text{$\cinfty$ or $\D'$}.
\end{equation}
In particular, $\HH_\xi^q(\cinfty(M)) = \HH_{\omega}^q(\cinfty(M))$ when $\omega=\xi \cdot \pmb{\omega}$.
If
$\omega$ is exact in $U$,
say
$\dd \phi = \omega|_U$ for a
$\phi \in \cinfty(U)$,
then
$\dd e^{i \phi} = i e^{i \phi} \omega|_U$
and
\begin{equation*}
  \dd (e^{i \phi} f)
  = \dd e^{i \phi} \wedge f + e^{i \phi} \dd f
  = e^{i \phi} \DD_\omega f
\end{equation*}
whatever
$f \in \D'(U; \Lambda^q)$.
In particular, for $q = 0$,
if $U$ is connected,
\begin{equation}
  \label{eq:solsDU}
  f \in \D'(U), \ \DD_\omega f = 0
  \Longrightarrow
  f = \mathrm{const.} \cdot e^{-i \phi},
\end{equation}
thus the sheaves of
homogeneous solutions of
$\DD_\omega$
in smooth functions
and in distributions
are one and the same;
we will denote this sheaf by
$\sol_\omega$.
A crucial feature is that
$\sol_\omega$
is \emph{locally} isomorphic with
the constant sheaf
$\sol_0$
-- its stalks are copies of $\C$ --
but not \emph{globally} in general
(see below).
We have commutative diagrams
\begin{equation*} 
  \begin{tikzcd}
    \D'(U; \Lambda^q) \arrow[r, "\DD_\omega"] \arrow[d, "e^{i \phi} \cdot"]
    &
    \D'(U; \Lambda^{q + 1}) \arrow[d, "e^{i \phi} \cdot"]  \\
    \D'(U; \Lambda^q) \arrow[r, "\dd"]
    & \D'(U; \Lambda^{q + 1}) 
  \end{tikzcd}
\end{equation*}
and the same holds
for smooth sections.
Since
$e^{i \phi} \neq 0$ on $U$,
we conclude that
\begin{equation}
  \label{eq:its-locally-deRham}
  \HH_\omega^q(\mathscr{F}(U)) \cong H_{\mathrm{dR}}^q(U),
  \quad \mathscr{F} = \text{$\cinfty$ or $\D'$},
\end{equation}
the right-hand side
standing for the usual
de Rham cohomology space.

Let
$\mathscr{F}_q$ denote the
sheaf of germs of
smooth or distributional sections
(depending on the case) of
$\Lambda^q$,
that is, forms and currents.
We have the exactness in degree
$q \in \{1, \ldots, n\}$ of the sequence
\begin{equation}
  \label{eq:sheaf-sequences-Domega}
  \begin{tikzcd}
    \cdots             \arrow[r, "\DD_\omega"] &
    \mathscr{F}_q      \arrow[r, "\DD_\omega"] &
    \mathscr{F}_{q + 1} \arrow[r, "\DD_\omega"] &
    \cdots .
  \end{tikzcd} 
\end{equation}
Indeed,
by the Poincar{\'e} Lemma,
every point
$t \in M$ has a neighborhood in which
$\omega$ is exact,
thus we can apply the
isomorphism~\eqref{eq:its-locally-deRham}
and use that the exterior derivative
$\dd$ is locally solvable at $t$ in degree $q \in \{1, \ldots, n\}$
(this is usually stated for smooth forms,
but also holds for currents).
In particular,
in both cases the
sequence~\eqref{eq:sheaf-sequences-Domega}
provides fine resolution of
$\sol_\omega$,
so on \emph{any} open set
$U$, we have
\begin{equation*}
  \HH_{\omega}^q(\cinfty(U))
  \cong \HH_{\omega}^q(\D'(U))
  \cong H^q(U; \sol_\omega),
\end{equation*}
the right-hand side denoting
the $q$-th cohomology space
with values in the sheaf
$\sol_\omega$.
Specifically,
the natural homomorphism of sheaves
$\Phi^0: \cinfty_0 \to \D'_0$
induces homomorphisms
$\Phi^q: \cinfty_q \to \D'_q$,
forming a homomorphism between resolutions
\begin{equation*}
  \begin{tikzcd}
    0              \arrow[r] &
    \sol_\omega     \arrow[r]               \arrow[d, "\mathrm{id}"] &
    \cinfty_0      \arrow[r, "\DD_\omega"]  \arrow[d, "\Phi^0"] &
    \cdots         \arrow[r, "\DD_\omega"]  &
    \cinfty_q      \arrow[r, "\DD_\omega"]  \arrow[d, "\Phi^q"] &
    \cinfty_{q + 1} \arrow[r, "\DD_\omega"]  \arrow[d, "\Phi^{q + 1}"] &
    \cdots \\
    0              \arrow[r] &
    \sol_\omega     \arrow[r] &
    \D'_0          \arrow[r, "\DD_\omega"] &
    \cdots         \arrow[r, "\DD_\omega"] &
    \D'_{q}        \arrow[r, "\DD_\omega"] &
    \D'_{q + 1}     \arrow[r, "\DD_\omega"] &
    \cdots
  \end{tikzcd}.
\end{equation*}
Hence, for each
$q \in \{0, \ldots, n\}$,
the map
\begin{equation}
  \label{eq:iso-sheaves-local}
  \Phi^q_*: \HH_{\omega}^q(\cinfty(U)) \lra \HH_{\omega}^q(\D'(U)),
  \quad U \sset M, 
\end{equation}
which is precisely
the one induced by the inclusion map
$\cinfty(U; \Lambda^q) \hookrightarrow \D'(U; \Lambda^q)$
on the quotients,
is an isomorphism\footnote{See
  e.g.~\cite[Theorem~3.13]{wells_dacm},
  but especially
  the remark by its end
  concerning naturality.}
-- we have just proved
a special case of
the so-called Atiyah-Bott Lemma.
Unwinding quotients
in~\eqref{eq:cohomology-spaces-Domega},
one deduces that:
\begin{Thm}
  \label{thm:analysis-fourier-cohomology}
  Given an open set
  $U \sset M$,
  we have that:
  \begin{enumerate}
  \item For every
    $f \in \cinfty(U; \Lambda^q)$
    such that there exists
    $u \in \D'(U; \Lambda^{q - 1})$
    with
    $\DD_\omega u = f$,
    there exists
    $v \in \cinfty(U; \Lambda^{q - 1})$
    such that
    $\DD_\omega v = f$.
    I.e.,
    $\DD_\omega|_U$ is~$\mathrm{(AGH)}$.
  \item For every
    $f \in \D'(U; \Lambda^q)$
    satisfying
    $\DD_\omega f = 0$,
    there exist
    $g \in \cinfty(U; \Lambda^q)$
    and
    $u \in \D'(U; \Lambda^{q - 1})$
    such that
    $f - g = \DD_\omega u$.
  \end{enumerate}
\end{Thm}

\begin{proof}
  Clearly,
  the first claim
  is equivalent to the injectivity
  of~\eqref{eq:iso-sheaves-local},
  while the second one
  is equivalent to its surjectivity --
  both of them established above.
\end{proof}

In particular,
$\DD_\omega$ is~$\mathrm{(AGH)}$
in each degree.
On time,
we recall that, by elliptic theory,
all the cohomology spaces
$\HH_{\omega}^q(\cinfty(M))$
are finite dimensional since $M$ is compact.
Both properties
will be used heavily, often without mention,
from here on.

Another consequence
of~\eqref{eq:solsDU}
is that, given a coordinate ball
$U \sset M$,
a function in
$\sol_\omega(U)$
either vanishes identically
or is never zero.
In particular,
the zero set of a global homogeneous solution
$f \in \sol_\omega(M)$
is both open and closed,
so by connectedness, such an $f$
also either vanishes everywhere
or not at all.
\begin{Lem}
  \label{lem:classification_sheaves}
  Either
  $\sol_\omega(M) = \{0\}$
  or
  $\dim_\C \sol_\omega(M) = 1$.
  The latter case happens
  if and only if
  $\sol_\omega$
  is isomorphic with
  the constant sheaf
  $\sol_0$.
\end{Lem}
\begin{proof}
  Suppose there exists
  $f \in \sol_\omega(M)$ non-zero.
  By the previous remarks,
  $f$ never vanishes.
  We will prove that any
  $g \in \sol_\omega(M)$
  is a constant times $f$.

  On an open set
  $U \sset M$ in which
  $\omega$ is exact,
  we have by~\eqref{eq:solsDU} that
  $f|_U = c_1 e^{-i \phi}$ and
  $g|_U = c_2 e^{-i \phi}$
  for some constants
  $c_1, c_2 \in \C$
  with $c_1 \neq 0$.
  Hence
  $(g / f)|_U = c_1^{-1} c_2$,
  and, since such open sets cover $M$,
  we reach the conclusion that
  $g / f$ is locally constant:
  by connectedness,
  this must be actually
  a constant function.

  Now,
  we prove that
  multiplication by $f$
  defines a sheaf isomorphism
  $\sol_0 \to \sol_\omega$,
  a claim that we check locally.
  Given $t \in M$,
  take a coordinate open ball $U \sset M$
  around it. We have
  \begin{equation*}
    h \in \sol_0(U)
    \Longrightarrow
    h = \mathrm{const.}
    \Longrightarrow
    h \cdot f|_U = \mathrm{const.} \cdot f|_U \in \sol_\omega(U),
  \end{equation*}
  so we have a monomorphism
  $\sol_0(U) \to \sol_\omega(U)$.
  But, since
  $\omega|_U$ has a primitive,
  any
  $g \in \sol_\omega(U)$
  is a multiple of
  $f|_U$,
  hence that monomorphism
  is also surjective.
\end{proof}

\begin{Cor}
  \label{cor:exacts_ZM}
  If
  $\omega$ is exact
  then
  $\sol_\omega \cong \sol_0$.
\end{Cor}
\begin{proof}
  If
  $\phi \in \cinfty(M)$
  is a primitive of
  $\omega$,
  then
  $\DD_{\omega} (e^{-i \phi}) = 0$,
  i.e.~$e^{-i \phi} \in \sol_\omega(M)$.
\end{proof}

\begin{Rem}
  The sheaf
  $\sol_\omega$
  is always locally isomorphic to
  $\sol_0$;
  the isomorphism can be made global
  precisely when
  $\sol_\omega$
  admits a non-vanishing global section.
\end{Rem}

\begin{Lem}
  \label{lem:monoid_structure_ZM}
  Given
  two closed $1$-forms
  $\omega_1, \omega_2$ on $M$
  such that
  $\sol_{\omega_1}(M) \neq \{0\} \neq \sol_{\omega_2}(M)$,
  we have that
  $\sol_{\omega_1 + \omega_2}(M) \neq 0$.
\end{Lem}
\begin{proof}
  If $f_1$ (resp.~$f_2$)
  is a section of
  $\sol_{\omega_1}$ (resp.~$\sol_{\omega_2}$)
  then
  $f_1 f_2$ is a section of
  $\sol_{\omega_1 + \omega_2}$.
  Indeed:
  \begin{equation*}
    \dd (f_1 f_2)
    = \dd f_1 \wedge f_2 + f_1 \wedge \dd f_2
    = -i  \omega_1 \wedge f_1 \wedge f_2
    - i f_1 \wedge \omega_2 \wedge f_2 
    = - i (\omega_1 + \omega_2) \wedge (f_1 f_2).
  \end{equation*}
  In particular, if
  $f_1 \in \sol_{\omega_1}(M)$
  and
  $f_2 \in \sol_{\omega_2}(M)$
  are both non-zero then
  $f_1 f_2 \in \sol_{\omega_1 + \omega_2}(M)$
  is non-zero.
\end{proof}

\begin{Thm}
  The set
  \begin{equation*}
    \mathcal{Z}_M
    \dfn
    \{ [\omega] \in H^1_{\mathrm{dR}}(M) \st \sol_\omega \cong \sol_0 \}
    =
    \{ [\omega] \in H^1_{\mathrm{dR}}(M) \st \sol_\omega(M) \neq \{0\} \}
  \end{equation*}
  is a group.
\end{Thm}
\begin{proof}
  Notice that
  $\mathcal{Z}_M$
  is well-defined:
  if
  $\omega, \omega^\bb$
  are two closed $1$-forms
  in the same cohomology class
  such that
  $\sol_\omega \cong \sol_0$
  then also
  $\sol_{\omega^\bb} \cong \sol_0$.
  Indeed,
  in that case there exists
  $\phi \in \cinfty(M)$
  such that
  $\omega^\bb = \omega + \dd \phi$:
  since
  $\sol_\omega \cong \sol_0$
  by assumption
  and
  $\sol_{\dd \phi} \cong \sol_0$
  by
  Corollary~\ref{cor:exacts_ZM},
  it follows from 
  Lemmas~\ref{lem:monoid_structure_ZM}
  and~\ref{lem:classification_sheaves}
  that
  $\sol_{\omega^\bb} \cong \sol_0$.

  The same argument
  shows that
  $\mathcal{Z}_M$
  is an additive submonoid of
  $H^1_{\mathrm{dR}}(M)$.
  Regarding inverses, if
  $[\omega] \in \mathcal{Z}_M$
  then there exists a non-vanishing
  $f \in \sol_\omega(M)$;
  clearly,
  \begin{equation*}
    \dd f + i f \omega = 0
    \Longrightarrow
    \dd (1/f) - i (1/f) \omega = 0.
  \end{equation*}
  Therefore,
  $1/f$
  is a non-vanishing element of
  $\sol_{-\omega}(M)$,
  so
  $-[\omega] \in \mathcal{Z}_M$.
\end{proof}

\begin{Rem}
  Recall that
  a real closed $1$-form
  $\alpha$ is \emph{integral}
  if
  $\int_\gamma \alpha \in 2 \pi \Z$
  for every $1$-cycle
  $\gamma$ in $M$.
  It follows
  from~\cite[Lemma~2.1]{bcp96}
  that
  \begin{equation}
    \label{eq:ZMbcp96}
    \mathcal{Z}_M
    = \{ [\omega] \in H^1_{\mathrm{dR}}(M) \st
    \text{$\Re \omega$ is integral and $\Im \omega$ is exact} \}.
  \end{equation}
\end{Rem}

\section{Formal characterization of
  the closure of $\ran \dd'$}
\label{sec:formal1}

In this section,
we address the issue of
formal solvability
(i.e.~at the level of
partial Fourier series)
and relate it with
the notion of
global solvability in degree $q\in \{1, \ldots, n\}$.
Now that we have
all the tools at our disposal,
the proofs are
pretty straightforward. 
\begin{Lem}
  \label{lem:formal_closure}
  For
  $f \in \cinfty(\Omega; \LLambda^{q})$,
  the following are equivalent:
  \begin{enumerate}
  \item $f$ belongs
    to the closure of
    $\ran \{\dd': \cinfty(\Omega; \LLambda^{q - 1}) \to \cinfty(\Omega; \LLambda^{q}) \}$;
  \item $\hat{f}_\xi$ is
    $\dd'_\xi$-exact
    for every $\xi \in \Z^m$.
  \end{enumerate}
\end{Lem}

\begin{Rem}
  The notion of
  $\hat{f}_\xi$
  being
  $\dd'_\xi$-exact
  is unambiguous:
  if
  $\dd'_\xi u_\xi = \hat{f}_\xi$
  is solvable
  in distributions,
  we can always find
  a smooth solution
  since
  $\dd'_\xi$ is~$\mathrm{(AGH)}$
  (Theorem~\ref{thm:analysis-fourier-cohomology}).
\end{Rem}
\begin{proof}
  Let
  $f \in \cinfty(\Omega; \LLambda^{q})$
  and suppose
  there exists
  $\{ u_\nu \}_{\nu \in \N} \sset \cinfty(\Omega; \LLambda^{q - 1})$
  such that
  $\dd' u_\nu \to f$
  in $\cinfty(\Omega; \LLambda^{q})$.
  By Proposition~\ref{prop:continuity_Fxi}
  we have that
  \begin{equation*}
    \dd'_\xi \mathcal{F_\xi} (u_\nu)
    =
    \mathcal{F_\xi} (\dd' u_\nu)
    \lra
    \hat{f}_\xi
    \quad \text{in $\cinfty(M; \Lambda^q)$}.
  \end{equation*}
  Therefore,
  $\hat{f}_\xi$
  belongs to the closure of
  $\ran \{ \dd'_\xi: \cinfty(M; \Lambda^{q - 1}) \to \cinfty(M; \Lambda^q) \}$,
  which is already closed
  in $\cinfty(M; \Lambda^q)$
  -- for instance,
  by Theorems~\ref{thm:agh_gs}
  and~\ref{thm:analysis-fourier-cohomology} --,
  hence:
  \begin{equation}
    \label{eq:formally_solvable}
    \hat{f}_\xi \in \ran \{\dd'_\xi: \cinfty(M; \Lambda^{q - 1}) \lra \cinfty(M; \Lambda^q) \},
    \quad \forall \xi \in \Z^m.
  \end{equation}
  Conversely, if
  $f \in \cinfty(\Omega; \LLambda^{q})$
  is such
  that~\eqref{eq:formally_solvable} holds,
  then
  for each $\xi \in \Z^m$
  there exists
  $u_\xi \in \cinfty(M; \Lambda^{q - 1})$
  such that
  $\dd'_\xi u_\xi = \hat{f}_\xi$,
  hence in the topology of
  $\cinfty(\Omega; \LLambda^{q})$
  we have
  \begin{equation*}
    f
    = \frac{1}{(2\pi)^m} \sum_{\xi \in \Z^m} e^{ix \xi} \wedge \dd'_\xi u_\xi
    = \lim_{\nu \to \infty} \dd' \Bigg( \frac{1}{(2\pi)^m} \sum_{|\xi| \leq \nu} e^{ix \xi} \wedge u_\xi \Bigg),
  \end{equation*}
  which proves that
  $f$ lies in the closure of
  $\ran \{ \dd': \cinfty(\Omega; \LLambda^{q - 1}) \to \cinfty(\Omega; \LLambda^{q}) \}$.
\end{proof}
\begin{Rem}
  \label{rem:compatibility_conditions}
  Lemma~\ref{lem:formal_closure}
  helps to settle the question
  of equivalence between
  different formulations of
  compatibility conditions appearing
  in the literature
  (hence, the equivalence between
  different notions of
  (smooth) global solvability).

  Consider, for instance,
  the space
  $E \sset \cinfty(\Omega; \LLambda^1)$
  determined by the compatibility conditions
  in~\cite{hz17, hz19}
  (which deal with the case
  $m = 1$, $q = 1$).
  Given
  $f \in \cinfty(\Omega; \LLambda^1)$,
  it follows easily from
  de Rham Theorem that,
  if
  $\hat{f}_\xi$ is
  $\dd'_\xi$-exact for every
  $\xi \in \Z^m$,
  then
  $f \in E$;
  whereas the converse
  follows from results
  in~\cite{hz17, hz19}:
  if
  $f \in E$,
  a solution to
  $\dd'_\xi u_\xi = \hat{f}_\xi$
  is obtained for every
  $\xi \in \Z^m$,
  first in some covering space of
  $M$,
  and then in
  $M$ by a convenient choice of
  initial conditions.
  In particular,
  $E$
  equals the closure of
  $\ran \{ \dd': \cinfty(\Omega) \to \cinfty(\Omega; \LLambda^{1}) \}$
  in that case.
\end{Rem}
The previous lemma yields 
our first major result.
\begin{Thm}
  \label{thm:main_agh}
  If~\eqref{eq:def_dprime}
  has closed range, then it
  is~$\mathrm{(AGH)}$.
\end{Thm}
\begin{proof}
  Suppose that
  $\dd': \cinfty(\Omega; \LLambda^{q}) \to \cinfty(\Omega; \LLambda^{q + 1})$
  has closed range
  and let
  $u \in \D'(\Omega; \LLambda^{q})$
  be such that
  $f \dfn \dd' u \in \cinfty(\Omega; \LLambda^{q + 1})$.
  Then
  $\hat{f}_\xi = \dd'_\xi \hat{u}_\xi$,
  that is,
  $\hat{f}_\xi$ is
  $\dd'_\xi$-exact
  for every $\xi \in \Z^m$.
  By Lemma~\ref{lem:formal_closure},
  $f$ belongs to
  $\ran \{ \dd': \cinfty(\Omega; \LLambda^{q}) \to \cinfty(\Omega; \LLambda^{q + 1}) \}$,
  so there exists
  $v \in \cinfty(\Omega; \LLambda^{q})$
  such that
  $\dd' v = f = \dd' u$.
\end{proof}
It generalizes~\cite[Corollary~7.2]{hz19}
(and preceding results):
\begin{Cor}
  \label{cor:main_agh}
  If the operator
  $\dd': \cinfty(\Omega) \to \cinfty(\Omega; \LLambda^1)$
  satisfies the property
  \begin{equation}
    \label{eq:trivial_kernel}
    \forall u \in \D'(\Omega), \ \dd' u = 0
    \Longrightarrow
    u = \mathrm{const.},
  \end{equation}
  then
  it has closed range
  if and only if
  it is globally hypoelliptic.
\end{Cor}
\begin{proof}
  This is an immediate consequence of
  Theorems~\ref{thm:agh_gs}
  and~\ref{thm:main_agh},
  since global hypoellipticity
  is equivalent
  to~$\mathrm{(AGH)}$
  when~\eqref{eq:trivial_kernel}
  holds.
\end{proof}
\begin{Rem}
  \label{rem:hz_case}
  Notice that~\eqref{eq:trivial_kernel}
  holds in the case considered
  in~\cite[Corollary~7.2]{hz19}
  ($m = 1$,
  $\omega = ib$ with
  $b$ real, closed and non-exact)
  thanks
  to~\cite[Lemma~2.2]{bcp96}.
  See also further discussion
  below about the group
  $\Gamma_{\pmb{\omega}}$.
\end{Rem}
\begin{Def}
  \label{def:formally_solvable}
  We say that an
  $f \in \cinfty(\Omega; \LLambda^{q})$
  is \emph{formally solvable}
  if
  for each $\xi \in \Z^m$
  there exists
  $u_\xi \in \cinfty(M; \Lambda^{q - 1})$
  such that
  $\dd'_\xi u_\xi = \hat{f}_\xi$.
\end{Def}
A sufficient condition
for that to happen is that
there exists
$u$ in
$\cinfty(\Omega; \LLambda^{q - 1})$
-- or even in
$\D'(\Omega; \LLambda^{q - 1})$ --
such that
$\dd' u = f$:
Lemma~\ref{lem:formal_closure}
entails the following converse.
\begin{Cor}
  \label{cor:formal_iff_closed}
  The following are equivalent:
  \begin{enumerate}
  \item $\dd': \cinfty(\Omega; \LLambda^{q - 1}) \to \cinfty(\Omega; \LLambda^{q})$
    has closed range;
  \item For every formally solvable
    $f \in \cinfty(\Omega; \LLambda^{q})$ 
    there exists
    $u \in \cinfty(\Omega; \LLambda^{q - 1})$
    such that
    $\dd' u = f$.
  \end{enumerate}
\end{Cor}

Before we move on,
we must introduce an
important set of frequencies.
For
$\pmb{\omega} = (\omega_1, \ldots, \omega_m)$,
the mapping
$\xi \in \Z^m \mapsto [\xi \cdot \pmb{\omega}] \in H^1_{\mathrm{dR}}(M)$
is a homomorphism of groups,
hence 
\begin{equation*}
  \Gamma_{\pmb{\omega}} \dfn \{ \xi \in \Z^m \st [\xi \cdot \pmb{\omega}] \in \mathcal{Z}_M \}
\end{equation*}
is a subgroup of $\Z^m$.
In particular,
either one of the
three mutually exclusive alternatives
hold:
\begin{enumerate}
\item $\Gamma_{\pmb{\omega}} = \{0 \}$;
\item $\Gamma_{\pmb{\omega}}$
  is infinite and proper;
  or
\item $\Gamma_{\pmb{\omega}} = \Z^m$.
\end{enumerate}
Using~\eqref{eq:ZMbcp96}
one has,
for instance,
that:
\begin{align*}
  \Gamma_{\pmb{\omega}} = \{0 \}
  &\Longleftrightarrow
  \text{$\Re (\xi \cdot \pmb{\omega})$ non-integral
    or $\Im (\xi \cdot \pmb{\omega})$ non-exact,
    $\forall \xi \in \Z^m \setminus \{0\}$
  } \\
  &\Longrightarrow
  \text{$\Re \omega_k$ non-integral
    or $\Im \omega_k$ non-exact,
    $\forall k \in \{1, \ldots, m\};$
  }
\end{align*}
while
\begin{align*}
  \Gamma_{\pmb{\omega}} = \Z^m
  &\Longleftrightarrow
  \text{$\Re (\xi \cdot \pmb{\omega})$ integral
    and $\Im (\xi \cdot \pmb{\omega})$ exact,
    $\forall \xi \in \Z^m$
  } \\
  &\Longleftrightarrow
  \text{$\Re \omega_k$ integral
    and $\Im \omega_k$ exact,
    $\forall k \in \{1, \ldots, m\}.$
  }
\end{align*}

\section{Applications to real structures}

Throughout this section,
$\omega_1, \ldots, \omega_m$
are assumed \emph{real}.
In this case,
\begin{equation*}
  \Gamma_{\pmb{\omega}}
  =
  \{ \xi \in \Z^m \st
  \text{$\xi \cdot \pmb{\omega}$ is integral} \}.
\end{equation*}

\subsection{Isomorphism theorems}
\label{sec:isomorphisms}

For
$X \subset \Z^m$
we set,
for each
$q \in \{0, \ldots, n \}$,
\begin{equation*}
  \cinfty_X (\Omega; \LLambda^{q})
  \dfn
  \{ f \in \cinfty (\Omega; \LLambda^{q} ) \st
  \hat f_\xi = 0, \ \forall \xi \in \Z^m \setminus X \},
\end{equation*}
and for
$f \in \cinfty(\Omega; \LLambda^{q})$,
we consider its projection on
$\cinfty_X (\Omega; \LLambda^{q})$:
\begin{equation*}
  f_X \dfn \frac{1}{(2\pi)^m} \sum_{\xi \in X} e^{ix \xi} \wedge \hat f_\xi.
\end{equation*}
It follows that
$f = f_X + f_{\Z^m \setminus X}$
and hence
\begin{equation*}
  \cinfty (\Omega; \LLambda^{q} )
  =
  \cinfty_X (\Omega; \LLambda^{q})
  \oplus
  \cinfty_{\Z^m \setminus X} (\Omega; \LLambda^{q}).
\end{equation*}
Notice moreover that
\begin{equation*}
  \dd' \cinfty_X (\Omega; \LLambda^{q}) \subset \cinfty_X (\Omega; \LLambda^{q + 1})
\end{equation*}
and that
\begin{equation*}
  \dd' u = f
  \Longleftrightarrow
  \text{
    $\dd' u_X = f_X$
    and
    $\dd' u_{\Z^m \setminus X} = f_{\Z^m \setminus X}$
    }.
\end{equation*}
It thus makes sense to define
\begin{equation*}
  \HH^{q}_{\dd'}(\cinfty_X(\Omega))
  \dfn
  \frac{ \ker \{ \dd': \cinfty_X(\Omega; \LLambda^{q}) \lra \cinfty_X(\Omega; \LLambda^{q + 1}) \} }
       { \ran \{ \dd': \cinfty_X(\Omega; \LLambda^{q - 1}) \lra \cinfty_X(\Omega; \LLambda^{q}) \} },
\end{equation*}
hence, we have
\begin{equation*}
  \HH^{q}_{\dd'}(\cinfty(\Omega))
  \cong
  \HH^{q}_{\dd'}(\cinfty_X(\Omega))
  \oplus
  \HH^{q}_{\dd'}(\cinfty_{\Z^m \setminus X}(\Omega)).
\end{equation*}

We denote by
$\pi: \tilde M \to M$
the universal covering of
$M$.
For each
$k \in \{1, \ldots, m\}$
there exists
$\psi_k \in \cinfty (\tilde M; \R)$
such that
$\dd \psi_k = \pi^* \omega_k$ on
$\tilde{M}$
(since the latter form is exact).
More generally, for 
$\xi \in \Z^m$
we set
\begin{equation}
  \label{definicaophixi}
  \psi_\xi \dfn \sum_{k = 1}^m \xi_k \psi_k \in \cinfty (\tilde M; \R),
\end{equation}
hence
$\dd \psi_\xi = \pi^* (\xi \cdot \pmb{\omega})$.
In an open set
$U \subset M$
in which
$\pi_U: \tilde{U} \to U$
is diffeomorphism
($\tilde{U} \sset \tilde M$
being another open set)
we have
\begin{equation*}
  \dd e^{i\psi_\xi \circ \pi_U^{-1}}
  =
  i e^{i\psi_\xi \circ \pi_U^{-1}} (\xi \cdot \pmb{\omega})
  \quad \text{on $U$}.
\end{equation*}
If
$\xi \in \Gamma_{\pmb{\omega}}$,
then $\xi \cdot \pmb{\omega}$ is integral,
hence by~\cite[Lemma~2.3]{bcm93},
if
$P,Q \in \tilde M$ satisfy
$\pi(P) = \pi(Q)$ then
$\psi_\xi(P) - \psi_\xi(Q) \in 2 \pi \Z$,
we can define for every
$\xi \in \Gamma_{\pmb{\omega}}$
a smooth function
\begin{equation}
  \label{estabemdefinidaemM}
  M \ni t \longmapsto e^{i \psi_\xi\circ \pi^{-1}(t)},
\end{equation}
even though
$\pi^{-1}$
is not a function.
It follows that
\begin{equation*}
   e^{i\psi_\xi\circ \pi^{-1}} \in \sol_{(- \xi \cdot \pmb{\omega})} (M) \setminus \{0\},
  \quad \forall \xi \in \Gamma_{\pmb{\omega}}.
\end{equation*}

We define a map
$\Psi: \tilde M \times \TT^m \to \tilde M \times \TT^m$
by
\begin{equation*}
  \Psi(\tilde t, x_1, \ldots, x_m)
  \dfn
  (\tilde t, x_1 - \psi_1(\tilde t), \ldots, x_m - \psi_m(\tilde t) ),
\end{equation*}
in which the sum
$x_k - \psi_k(\tilde{t})$
takes place in
$\R/(2\pi \Z)$.
In particular,
$\Psi$ is a diffeomorphism with inverse
\begin{equation*}
  \Psi^{-1} (\tilde t, x_1, \ldots, x_m )
  =
  (\tilde t, x_1 + \psi_1(\tilde t), \ldots, x_m + \psi_m(\tilde t) ).
\end{equation*}
\begin{Prop}
  \label{Prop:freq-desc}
  Suppose that
  $f \in \cinfty_{\Gamma_{\pmb{\omega}}} (\Omega; \LLambda^{q})$.
  Then
  \begin{equation}
    \label{eq:def_Theta}
    \Theta (f)
    \dfn
    \frac{1}{(2\pi)^m}
    \sum_{\xi \in \Gamma_{\pmb{\omega}}} e^{ix \xi- i\psi_\xi\circ \pi^{-1}} \hat f_\xi
  \end{equation}
  defines an element in
  $\cinfty_{\Gamma_{\pmb{\omega}}} (\Omega; \LLambda^{q})$.
  The map
  \begin{equation*}
    \Theta:
    \cinfty_{\Gamma_{\pmb{\omega}}} (\Omega; \LLambda^{q})
    \lra
    \cinfty_{\Gamma_{\pmb{\omega}}} (\Omega; \LLambda^{q})
  \end{equation*}
  is a
  linear isomorphism
  and satisfies
  $\dd' \circ \Theta = \dd_t \circ \Theta$.
\end{Prop}
\begin{proof}
  Writing the partial Fourier series of
  $f$, we have that
  \begin{equation*}
    f(t,x)
    =
    \frac{1}{(2\pi)^m} 
    \sum_{\xi \in \Gamma_{\pmb{\omega}}} e^{ix \xi} \hat{ f}_\xi(t).
  \end{equation*}
  We define
  $F \dfn (\pi \times \mathrm{id}_{\TT^m})^* f$,
  which is a smooth
  $q$-form on
  $\tilde{M} \times \TT^m$.
  By continuity of the pullback map,
  it follows that
  \begin{equation*}
    F(\tilde{t}, x)
    =
    \frac{1}{(2\pi)^m}
    \sum_{\xi \in \Gamma_{\pmb{\omega}}} e^{ix \xi} (\pi^* \hat{f}_\xi)(\tilde{t}),
  \end{equation*}
  with convergence in
  the space of smooth
  $q$-forms on
  $\tilde{M} \times \TT^m$.
  If we change coordinates
  using the diffeomorphism
  $\Psi$, we obtain
  \begin{equation*}
    (\Psi^* F)(\tilde{t}, x)
    =
    \frac{1}{(2\pi)^m}
    \sum_{\xi \in \Gamma_{\pmb{\omega}}} e^{i(x \xi- \psi_\xi(\tilde{t}))} (\pi^* \hat{f}_\xi)(\tilde{t}),
  \end{equation*}
  which we will show
  that descends to
  $M \times \TT^m$ as
  $\Theta(f)$.
  The sequence of
  truncated sums
  \begin{equation*}
    f_\nu(t)
    \dfn
    \frac{1}{(2\pi)^m}
    \sum_{\substack{\xi \in \Gamma_{\pmb{\omega}} \\ |\xi| \leq \nu}} e^{ix \xi} \hat{f}_\xi (t),
    \quad \nu \in \N,
  \end{equation*}
  which approximates
  $f$ in
  $\cinfty(\Omega; \LLambda^q)$,
  certainly satisfy
  \begin{equation*}
    \Psi^* (\pi \times \mathrm{id}_{\TT^m})^* f_\nu
    \lra
    \Psi^* F
    \quad
    \text{on $\tilde{M} \times \TT^m$}
  \end{equation*}
  whereas
  \begin{equation*}
    g_\nu(t)
    \dfn
    \frac{1}{(2\pi)^m}
    \sum_{\substack{\xi \in \Gamma_{\pmb{\omega}} \\ |\xi| \leq \nu}} e^{ix \xi- i\psi_\xi\circ \pi^{-1}(t)} \hat f_\xi(t),
    \quad \nu \in \N,
  \end{equation*}
  satisfy
  $(\pi \times \mathrm{id}_{\TT^m})^* g_\nu = \Psi^* (\pi \times \mathrm{id}_{\TT^m})^* f_\nu$
  clearly,
  hence
  \begin{equation}
    \label{eq:a_convergence}
    (\pi \times \mathrm{id}_{\TT^m})^* g_\nu
    \lra
    \Psi^* F
    \quad
    \text{on $\tilde{M} \times \TT^m$}.
  \end{equation}
  To prove convergence of
  $\{ g_\nu \}_{\nu \in \N}$ in
  $\cinfty(\Omega; \LLambda^q)$,
  it suffices to check it on
  $U \times \TT^m$
  for a suitably small open set $U \sset M$:
  its limit will then be automatically
  a globally defined
  $q$-form
  $\Theta(f)$
  on
  $M \times \TT^m$,
  and can be easily seen to
  be a section of
  $\LLambda^q$.
  This is the case if,
  say,
  there exists an open set
  $\tilde{U} \sset \tilde{M}$
  such that
  $\pi_U: \tilde{U} \to U$
  is a diffeomorphism. 
  Since the previous
  convergence~\eqref{eq:a_convergence}
  takes place in
  $\tilde{U} \times \TT^m$
  as well by restriction,
  in which
  $\pi_U \times \mathrm{id}_{\TT^m}$
  is invertible,
  we conclude that
  $\{ g_\nu \}_{\nu \in \N}$
  converges in
  $U \times \TT^m$,
  proving our claim.

  Its clear
  that~\eqref{eq:def_Theta} holds
  (as
  $\Theta(f)$
  is by definition the limit of
  $\{ g_\nu \}_{\nu \in \N}$),
  hence of course
  $\Theta (f) \in \cinfty_{\Gamma_{\pmb{\omega}}} (\Omega; \LLambda^{q})$:
  by continuity of each map
  $\mathcal{F}_\xi$,
  equation~\eqref{eq:def_Theta}
  gives its Fourier coefficients explicitly.
  Moreover:
  \begin{align*}
    (2\pi)^m \dd' \Theta (f)  &=
    \bigg( \dd_t + \sum_{k = 1}^m \omega_k \wedge \frac{\del}{\del x_k} \bigg)
    \sum_{\xi \in \Gamma_{\pmb{\omega}}} e^{ix \xi -i\psi_\xi\circ \pi^{-1}} \hat f_\xi  \\
    &=
    \sum_{\xi \in \Gamma_{\pmb{\omega}}}  \bigg(  e^{ix \xi}  \dd_t( e^{ -i\psi_\xi\circ \pi^{-1}} \hat{ f_\xi}) 
    +
    \sum_{k = 1}^m   \frac{\del e^{ix \xi}}{\del x_k} \omega_k\wedge ( e^{-i\psi_\xi\circ \pi^{-1}} \hat f_\xi )  \bigg) \\
    &=
    \sum_{\xi \in \Gamma_{\pmb{\omega}}}  \bigg(  e^{ix \xi-i \psi_\xi\circ \pi^{-1}} \Big( \dd_t \hat f_\xi - i (\xi \cdot \pmb{\omega}) \wedge \hat f_\xi \Big)   +  i (\xi \cdot \pmb{\omega}) \wedge e^{ix \xi- i\psi_\xi\circ \pi^{-1}} \hat f_\xi   \bigg) \\
    &=
    \sum_{\xi \in \Gamma_{\pmb{\omega}}}  e^{ix \xi- i\psi_\xi\circ \pi^{-1}} \widehat{\dd_t f_\xi}  \\
    &= (2\pi)^m \Theta (\dd_t f),
  \end{align*}
  that is:
  $\Theta^{-1} \circ \dd' \circ \Theta = \dd_t$
  on
  $\cinfty_{\Gamma_{\pmb{\omega}}} (\Omega; \LLambda^{q})$.
\end{proof}
A similar calculation shows that, if 
\begin{equation*}
  \Theta_\xi :  \cinfty (M; \Lambda^q)   \lra  \cinfty (M; \Lambda^q),  \quad
  \Theta_\xi (h)  \dfn  e^{i\psi_\xi\circ \pi^{-1}} h,
\end{equation*}
then
$\Theta_\xi^{-1} \circ \dd'_\xi \circ \Theta_\xi = \dd$,
whatever
$\xi \in \Gamma_{\pmb{\omega}}$.

\begin{Thm}
  \label{Thm:reduction-glob-solv}
  The map
  $\dd' : \cinfty (\Omega; \LLambda^{q}) \to \cinfty (\Omega; \LLambda^{q + 1})$
  has closed range
  if and only if
  given a formally solvable
  $f \in \cinfty_{\Z^m \setminus \Gamma_{\pmb{\omega}}} (\Omega; \LLambda^{q + 1})$,
  there exists
  $u \in \cinfty(\Omega; \LLambda^{q})$
  such that
  $\dd'u = f$.
\end{Thm}
\begin{proof}
  The direct implication
  is granted
  by Corollary~\ref{cor:formal_iff_closed};
  we prove the converse.
  Let
  $f \in \cinfty (\Omega; \LLambda^{q + 1})$
  be such that for each
  $\xi \in \Z^m$
  there exists
  $u_\xi \in \cinfty(M; \Lambda^q)$
  satisfying
  $\dd'_\xi u_\xi = \hat f_\xi$.
  Write
  $f = f_{\Gamma_{\pmb{\omega}}} + f_{\Z^m \setminus \Gamma_{\pmb{\omega}}}$
  and notice that both
  $f_{\Gamma_{\pmb{\omega}}}$ and
  $f_{\Z^m \setminus \Gamma_{\pmb{\omega}}}$
  are also formally solvable.
  By hypothesis,
  there exists
  $v \in \cinfty(\Omega; \LLambda^{q})$
  such that
  $\dd' v = f_{\Z^m \setminus \Gamma_{\pmb{\omega}}}$:
  we must solve
  $\dd' u = f_{\Gamma_{\pmb{\omega}}}$.
  
  We claim that
  $g \dfn \Theta^{-1} f_{\Gamma_{\pmb{\omega}}} \in \cinfty_{\Gamma_{\pmb{\omega}}} (\Omega; \LLambda^{q + 1} )$
  is formally solvable
  w.r.t.~$\dd_t$ in
  $\Omega$.
  Indeed,
  $\hat g_\xi = 0$
  for
  $\xi \in \Z^m \setminus \Gamma_{\pmb{\omega}}$,
  while for
  $\xi \in \Gamma_{\pmb{\omega}}$,
  we have
  \begin{equation*}
    \dd'_\xi u_\xi = \hat{f}_\xi
    \Longrightarrow
    \dd_t \Theta^{-1}_\xi u_\xi
    =
    ( \Theta_\xi^{-1} \dd_\xi' \Theta_\xi) \Theta^{-1}_\xi u_\xi
    =
    \Theta_\xi^{-1} \hat f_\xi
    =
    \mathcal{F}_\xi(\Theta^{-1} f)
    =
    \hat g_\xi.
  \end{equation*}
  But global solvability of $\dd_t$ is a general fact -- see Lemma~\ref{lem:GS_dt} below --  which will also play a role  later on.
  It then follows from
  Corollary~\ref{cor:formal_iff_closed} applied to
  $\dd_t$
  that there exists
  $w \in \cinfty(\Omega; \LLambda^{q})$
  such that
  $\dd_t w = g$.
  Note that we can replace
  $w$ with
  $w_{\Gamma_{\pmb{\omega}}}$,
  if necessary, and assume that
  $w \in \cinfty_{\Gamma_{\pmb{\omega}}}(\Omega; \LLambda^{q})$.
  If we set
  $u \dfn \Theta w \in \cinfty_{\Gamma_{\pmb{\omega}}}(\Omega; \LLambda^{q})$
  then
  \begin{equation*}
    \dd' u
    =
    \dd' \Theta w
    =
    \Theta \dd_t w
    = \Theta g
    = f_{\Gamma_{\pmb{\omega}}}
  \end{equation*}
  and we are done.
\end{proof}

\begin{Lem}
  \label{lem:GS_dt}
  $\dd_t : \cinfty(\Omega; \LLambda^{q}) \to \cinfty(\Omega; \LLambda^{q + 1})$
  has closed range.
\end{Lem}
\begin{proof}
  We endow
  $M$ with a Riemannian metric
  and let
  $\Delta \dfn \dd \dd^* + \dd^* \dd$
  be the Laplace-Beltrami operator
  acting on forms on $M$;
  it is elliptic of order $2$
  and therefore satisfies
  elliptic estimates:
  given
  $k \in \Z_+$
  there exists
  $C_k > 0$
  such that
  \begin{equation*}
    \| \psi \|_{\sob^{k+2} (M; \Lambda^q)}
    \leq
    C_k \left( \| \Delta \psi \|_{\sob^{k} (M; \Lambda^q)} + \| \psi \|_{\sob^{k} (M; \Lambda^q)} \right),
    \quad \forall \psi \in \sob^{k+2} (M; \Lambda^q),
  \end{equation*}
  where
  $\sob^k$ are Sobolev spaces.
  A standard argument
  shows that
  \begin{equation*}
    \| \psi \|_{\sob^{k+2} (M; \Lambda^q)}
    \leq
    C_k \| \Delta \psi \|_{\sob^{k} (M; \Lambda^q)},
    \quad \text{$\forall \psi \in \sob^{k+2} (M; \Lambda^q)$, $L^2$-orthogonal to $\ker \Delta$}.
  \end{equation*}
  The following orthogonal decomposition
  w.r.t.~the $L^2(M)$ metric
  is also well-known:
  \begin{equation}
    \label{eq:orth_decomp}
    \cinfty (M; \Lambda^q)
    =
    \ker \Delta \oplus \ran \dd \oplus \ran \dd^*
  \end{equation}
  (all operators acting on smooth forms),
  hence in particular,
  for
  $\psi \in \ran \dd^*$,
  \begin{equation*}
    \| \psi \|_{\sob^{k + 2} (M; \Lambda^q)}
    \leq
    C_k \| \dd^* \dd \psi \|_{\sob^{k} (M; \Lambda^q)}
    \leq
    C_k' \| \dd \psi \|_{\sob^{k - 1} (M; \Lambda^{q + 1})}
  \end{equation*}
  since in that case
  $\Delta \psi = \dd^* \dd \psi$.
  We conclude that,
  for each
  $k \in \Z_+$,
  there exist
  $c_k > 0$ and
  $j \in \Z_+$
  such that
  \begin{equation}
    \label{eq:estimate_closed_range2}
    \| \psi \|_{\sob^{k} (M; \Lambda^q)}
    \leq
    c_k \| \dd \psi \|_{\sob^{j} (M; \Lambda^{q + 1})},
    \quad \forall \psi \in \ran \dd^*.
  \end{equation}
  
  Now, we obtain our conclusion
  by means of 
  Theorem~\ref{thm:agh_gs}.
  Let
  $u \in \D'(\Omega; \LLambda^{q})$
  be such that
  $\dd_t u \in \cinfty(\Omega; \LLambda^{q + 1})$.
  Then,
  for every
  $\xi \in \Z^m$,
  we have that
  $\dd \hat{u}_\xi \in \cinfty(M; \Lambda^{q + 1})$
  by
  Corollary~\ref{cor:actions_fouriercoeff},
  hence, by
  Theorem~\ref{thm:analysis-fourier-cohomology},
  there exists
  $v_\xi \in \cinfty(M; \Lambda^q)$
  such that
  $\dd v_\xi = \dd \hat{u}_\xi$;
  thanks to~\eqref{eq:orth_decomp}
  we can further assume that
  $v_\xi \in \ran \dd^*$.
  By~\eqref{eq:estimate_closed_range2},
  for each
  $k \in \Z_+$,
  there exist
  $c_k > 0$ and
  $j \in \Z_+$
  such that
  \begin{equation*}
    \| v_\xi \|_{\sob^{k}(M; \Lambda^q)}
    \leq
    c_k \| \dd v_\xi \|_{\sob^j(M; \Lambda^{q + 1})},
    \quad \forall \xi \in \Z^m. 
  \end{equation*}
  Moreover, since
  $\dd_t u$ is smooth
  and
  $\mathcal{F}_\xi(\dd_t u) = \dd \hat{u}_\xi$,
  for every
  $j \in \Z_+$
  and
  $s > 0$,
  there exists a constant
  $A_{j,s} > 0$ such that
  \begin{equation*}
    \| \dd \hat{u}_\xi \|_{\sob^j(M; \Lambda^{q + 1})}
    \leq
    A_{j,s} (1 + |\xi|^2)^{-s},
    \quad \forall \xi \in \Z^m,
  \end{equation*}
  and so,
  for
  $c'_{k,s} \dfn c_k A_{j,s}$, we have
  \begin{equation}
    \label{eq:vxi_estimate}
    \| v_\xi \|_{\sob^{k}(M; \Lambda^{q})}
    \leq
    c'_{k,s} (1 + |\xi|^2)^{-s},
    \quad \forall \xi \in \Z^m.
  \end{equation}
  The latter ensures that
  the series
  \begin{equation*}
    v \dfn \frac{1}{(2 \pi)^m} \sum_{\xi \in \Z^m} e^{ix \xi} \wedge v_\xi
  \end{equation*}
  converges in
  $L^2(\Omega; \LLambda^{q})$
  since
  \begin{align*}
    \frac{1}{(2 \pi)^m}
    \sum_{\xi} \int_\Omega \| e^{i x \xi} v_\xi(t) \|_{\Lambda^q_t}^2 \dd V(t,x)
    &= \sum_{\xi} \int_M \| v_\xi(t) \|_{\Lambda^q_t}^2 \dd V_M(t) \\
    &\leq (c'_{0, 2m})^2 \sum_{\xi} (1 + |\xi|^2)^{-2m}
    < \infty.
  \end{align*}
  Moreover, for each
  $\xi \in \Z^m$, we have that
  $\hat{v}_\xi = v_\xi$,
  hence
  $\dd \hat{v}_\xi = \dd \hat{u}_\xi$,
  so
  $\dd_t v = \dd_t u$.
  Estimates~\eqref{eq:vxi_estimate}
  further ensure that
  $v \in \cinfty(\Omega; \LLambda^{q})$.
\end{proof}

\begin{Rem}
  We have proved, in
  Theorem~\ref{Thm:reduction-glob-solv},
  that
  $\dd'$ is always
  \emph{$\Gamma_{\pmb{\omega}}$-globally solvable},
  that is, if
  $f \in \cinfty_{\Gamma_{\pmb{\omega}}} (\Omega; \LLambda^{q + 1})$
  is formally solvable,
  then there exists
  $u \in \cinfty(\Omega; \LLambda^{q})$
  such that
  $\dd'u = f$,
  with no further hypotheses.
\end{Rem}

\subsubsection{Reduction in cohomology}
\label{sec:reduction}

The map
$\Theta$ descends to a linear isomorphism
\begin{equation*}
  \HH^{q}_{\dd'}(\cinfty_{\Gamma_{\pmb{\omega}}}(\Omega))
  \cong
  \HH^{q}_{\dd_t}(\cinfty_{\Gamma_{\pmb{\omega}}}(\Omega))
\end{equation*}
thanks to its properties
deduced above.
Now,
we provide a more detailed description
of those spaces,
for which we introduce some notation.
Take
$\xi^{(1)}, \ldots, \xi^{(r)}$
a basis of
$\Gamma_{\pmb{\omega}}$
as a
$\Z$-module\footnote{One should
pay special attention
to the case
$r = 0$,
i.e.~$\Gamma_{\pmb{\omega}} = \{0\}$.
When properly interpreted
(e.g.~$\TT^0 = \{1\}$,
$\cinfty(\TT^0) = \C$, etc.),
most of the results below
are trivial in that case.},
thus, for each
$\xi \in \Gamma_{\pmb{\omega}}$,
there exists a unique
$\eta \in \Z^{r}$
such that
\begin{equation*}
  \xi
  =
  \eta_1 \xi^{(1)} + \cdots + \eta_r \xi^{(r)}
  \dfn
  \eta \cdot \xi^{0}
\end{equation*}
in which
$\xi^{0} \dfn (\xi^{(1)}, \ldots, \xi^{(r)})$. We also define a smooth map $\theta: \TT^m \to \TT^r$
by
\begin{equation*}
  \theta(x_1, \ldots, x_m)   \dfn   \left ( x \cdot \xi^{(1)}, \ldots,  x \cdot\xi^{(r)} \right).
\end{equation*}
Note that
\begin{align*}
  \sum_{j=1}^{r} (x \cdot \xi^{(j)}) \eta_j = x_1\Big( \sum_{j=1}^{r} \xi_1^{(j)} \eta_j\Big) + \cdots + x_m \Big( \sum_{j=1}^{r} \xi_m^{(j)} \eta_j \Big)= x \cdot \xi
\end{align*}
thus
\begin{equation*}
  e^{i \theta(x) \cdot \eta}  =  \prod_{j = 1}^r (e^{i x \cdot \xi^{(j)}} )^{\eta_j}
  =  e^{i x \cdot (\eta \cdot \xi^0)}= e^{i x \xi}.
\end{equation*}
\begin{Prop}
  \label{prop:pullback_Tr}
  The pullback
  $\theta^*: \cinfty(\TT^r) \to \cinfty_{\Gamma_{\pmb{\omega}}}(\TT^m)$
  is a linear isomorphism.
\end{Prop}
\begin{proof}
  We write an
  $f \in \cinfty(\TT^r)$
  as
  \begin{equation*}
    f
    = 
    \frac{1}{(2 \pi)^r}
    \sum_{\eta \in \Z^r} e^{i y \eta} \hat{f}_{\eta}
  \end{equation*}
  with convergence in
  $\cinfty(\TT^r)$.
  Hence
  \begin{equation*}
    \theta^* f
    =
    \frac{1}{(2 \pi)^r}
    \sum_{\eta \in \Z^r} e^{i \theta(x) \cdot \eta} \hat{f}_{\eta}
    =
    \frac{1}{(2 \pi)^r}
    \sum_{\eta \in \Z^r} e^{i x (\eta \cdot \xi^0)} \hat{f}_{\eta}
  \end{equation*}
  meaning that
  \begin{equation*}
    \mathcal{F}_\xi (\theta^* f)
    =
    \begin{cases}
      (2 \pi)^{m - r} \hat{f}_{\eta},
      &\text{if $\xi = \eta \cdot \xi^0 \in \Gamma_{\pmb{\omega}}$}; \\
      0,
      &\text{if $\xi \notin \Gamma_{\pmb{\omega}}$}.
    \end{cases}
  \end{equation*}
  In particular,
  $\theta^* f \in \cinfty_{\Gamma_{\pmb{\omega}}}(\TT^m)$.
  Moreover,
  $\theta^*$ is injective since
  $\theta^* f = 0$ implies that
  $\hat{f}_{\eta} = 0$ for all
  $\eta \in \Z^r$.
  As for surjectivity,
  given
  $g \in \cinfty_{\Gamma_{\pmb{\omega}}}(\TT^m)$,
  let
  \begin{equation*}
    f
    \dfn 
    \frac{1}{(2 \pi)^r} \sum_{\eta \in \Z^r} e^{i y \eta} \hat{g}_{\eta \cdot \xi^0}
    \in \cinfty(\TT^r)
  \end{equation*}
  which clearly satisfies
  $\theta^* f = g$.

  Note that
  $|\xi|$ and $|\eta|$
  are comparable,
  so a series converges in
  $\cinfty(\TT^r)$ if and only if its image by
  $\theta^*$ converges in
  $\cinfty_{\Gamma_{\pmb{\omega}}}(\TT^m)$.
  Indeed,
  we can complete
  $\xi^{0}$ to a basis of $\R^m$,
  so we can write
  $A = (\xi^{(1)}, \ldots, \xi^{(r)}, \xi^{(r+1)}, \ldots, \xi^{(m)})$,
  in which
  $A$ is a matrix whose
  $j$-th column equals
  $\xi^{(j)}$.
  It follows that
  multiplication of
  $A$ by
  $(\eta, 0) \in \Z^{m}$ yields
  \begin{equation*}
    (\eta, 0) \cdot A = \eta \cdot \xi^{0} = \xi.
  \end{equation*}
  Since
  $A$ is invertible,
  it follows that
  $c |\eta| \leq |\xi| \leq C |\eta|$
  for some constants
  $c, C > 0$.
\end{proof}

Now given
$q \in \{0, \ldots, n\}$,
for any
$f \in \cinfty(\Omega; \LLambda^{q})$
and
$x \in \TT^m$,
we define
$f(x) \in \cinfty(M; \Lambda^q)$
by ``fixing the $x$-variable'',
i.e.~$f(t,x) = f(x)(t)$
for every
$t \in M$.
As such,
if
$\dd_t f = 0$,
then
$\dd f(x) = 0$,
and, if we pick
closed forms
$\tau_1, \ldots, \tau_{b_q} \in \cinfty(M; \Lambda^q)$
such that
$[\tau_1], \ldots, [\tau_{b_q}]$
is a basis of
$\HH^q_{\mathrm{dR}}(M)$,
then for each
$x \in \TT^m$
we write
\begin{equation*}
  [f(x)] = \sum_{\ell = 1}^{b_q} a_\ell(x) [\tau_\ell]
  \quad \text{in $\HH^q_{\mathrm{dR}}(M)$}
\end{equation*}
for some uniquely determined coefficients
$a_1(x), \ldots, a_{b_q}(x)$.
Things can be arranged so that
$x \in \TT^m \mapsto a_\ell(x) \in \C$
are all smooth.

Indeed,
by endowing
$M$ with a Riemannian metric,
we may pick
$\tau_1, \ldots, \tau_{b_q}$
as a basis of
$\ker \Delta$,
the space of harmonic
$q$-forms on
$M$,
which is orthonormal
w.r.t.~the $L^2$ inner product on
$\cinfty(M; \Lambda^q)$.
Then there exists
$u_x \in \cinfty(M; \Lambda^{q - 1})$
such that
\begin{equation*}
  f(x)
  =
  \sum_{\ell = 1}^{b_q} a_\ell(x) \tau_\ell
  +
  \dd_t u_x,
\end{equation*}
thus realizing the last sum
as the orthogonal projection of
$f(x)$ onto
$\ker \Delta$
(see~\eqref{eq:orth_decomp}):
\begin{equation}
  \label{eq:a_ell}
  a_\ell(x)
  =
  \langle f(x), \tau_\ell \rangle_{L^2(M; \Lambda^q)}
  =
  \int_M \langle f(t,x), \tau_\ell(t) \rangle_{\Lambda^q_t} \ \dd V_M(t),
  \quad \ell \in \{1, \ldots, b_q\},
\end{equation}
are therefore smooth
w.r.t.~$x$.

Notice also that
$a_1(x), \ldots, a_{b_q}(x)$
depend only on the cohomology class of
$f$:
if
$f^\bb \dfn f + \dd_t v$
for some
$v \in \cinfty(\Omega; \LLambda^{q - 1})$
then
\begin{equation*}
  f^\bb(x)
  =
  f(x) + \dd_t v(x)
  =
  \sum_{\ell = 1}^{b_q} a_\ell(x) \tau_\ell
  +
  \dd_t (u_x + v(x)).
\end{equation*}
It follows
from~\eqref{eq:a_ell}
that, for every
$\xi \in \Z^m$,
\begin{equation}
  \label{eq:a_ell_fourier}
  \mathcal{F}_\xi(a_\ell)
  =
  \int_{\TT^m} e^{-ix \xi} \int_M \langle f(t,x), \tau_\ell(t) \rangle_{\Lambda^q_t} \ \dd V_M(t) \dd x
  =
  \langle \hat{f}_\xi, \tau_\ell \rangle_{L^2(M; \Lambda^q)}.
\end{equation}

If we further assume that
$f \in \cinfty_{\Gamma_{\pmb{\omega}}} (\Omega; \LLambda^{q})$,
then by~\eqref{eq:a_ell_fourier}
we have
$\mathcal{F}_\xi(a_\ell) = 0$
whenever
$\xi \notin \Gamma_{\pmb{\omega}}$,
that is,
$a_\ell \in \cinfty_{\Gamma_{\pmb{\omega}}} (\TT^m)$,
so, by
Proposition~\ref{prop:pullback_Tr},
there is a unique
$a_\ell^\bb \in \cinfty(\TT^r)$
such that
$a_\ell = \theta^* (a_\ell^\bb)$.
We define
\begin{equation*}
  \mathsf{T}:
  \HH^{q}_{\dd_t}(\cinfty_{\Gamma_{\pmb{\omega}}}(\Omega))
  \lra
  \cinfty(\TT^r) \otimes \HH^q_{\mathrm{dR}}(M)
\end{equation*}
by
\begin{equation*}
  \mathsf{T}([f])
  \dfn \sum_{\ell = 1}^{b_q} a_\ell^\bb \otimes [\tau_\ell].
\end{equation*}
\begin{Thm}
  The map $\mathsf{T}$
  is a linear isomorphism.
\end{Thm}
\begin{proof}
  The surjectivity is obvious;
  given
  $a_1^\bb, \ldots, a_{b_q}^\bb \in \cinfty(\TT^r)$,
  we have
  \begin{equation*}
    f \dfn \sum_{\ell = 1}^{b_q} (\theta^* a_\ell^\bb) \tau_\ell
    \Longrightarrow
    \mathsf{T}([f]) = \sum_{\ell = 1}^{b_q} a_\ell^\bb \otimes [\tau_\ell].
  \end{equation*}
  
  As for injectivity,
  if a
  $\dd_t$-closed
  $f \in \cinfty_{\Gamma_{\pmb{\omega}}}(\Omega; \LLambda^{q})$
  is such that
  $\mathsf{T}([f]) = 0$, then
  $a_1 = \cdots = a_{b_q} = 0$.
  Since
  $\dd_t f = 0$,
  it follows that
  $\dd \hat{f}_\xi = 0$
  for every
  $\xi \in \Gamma_{\pmb{\omega}}$,
  and so we can find
  $a^\xi_\ell \in \C$ and
  $u_{\xi}\in \cinfty(M; \Lambda^{q - 1})$
  such that
  \begin{equation*}
    \hat{f}_\xi = \sum_{\ell = 1}^{b_q} a^{\xi}_\ell \tau_\ell + \dd u_{\xi},
    \quad \forall \xi \in \Gamma_{\pmb{\omega}}.
  \end{equation*}
  Note that for
  $\xi \notin \Gamma_{\pmb{\omega}}$
  we have
  $\hat{f}_\xi = 0$.
  By~\eqref{eq:a_ell_fourier},
  we have,
  for each
  $\ell \in \{1, \ldots, b_q\}$
  and
  $\xi \in \Gamma_{\pmb{\omega}}$,
  using that
  $\tau_\ell \in \ker \Delta = \ker \dd \cap \ker \dd^*$:
  \begin{equation*}
    \mathcal{F}_\xi (a_\ell)
    =
    a^{\xi}_\ell + \langle \dd u_{\xi}, \tau_\ell \rangle_{L^2(M; \Lambda^q)}
    =
    a^{\xi}_\ell.
  \end{equation*}
  By the assumption that
  $a_\ell = 0$ for every
  $\ell \in \{1, \ldots, b_q\}$,
  we reach the conclusion that
  $\hat{f}_\xi$ is exact for every
  $\xi \in \Z^m$,
  that is,
  $f$ is formally
  $\dd_t$-solvable.
  By
  Corollary~\ref{cor:formal_iff_closed}
  and
  Lemma~\ref{lem:GS_dt},
  we conclude that
  $f$ is
  $\dd_t$-exact,
  i.e.~$[f] = 0$
  in
  $\HH^{q}_{\dd_t}(\cinfty_{\Gamma_{\pmb{\omega}}}(\Omega))$.
\end{proof}

\subsubsection{Summary of the section}

We always have linear isomorphisms
\begin{equation*}
  \HH^{q}_{\dd'}(\cinfty(\Omega))
  \cong
  \HH^{q}_{\dd'}(\cinfty_{\Gamma_{\pmb{\omega}}}(\Omega))
  \oplus
  \HH^{q}_{\dd'}(\cinfty_{\Z^m \setminus \Gamma_{\pmb{\omega}}}(\Omega))
\end{equation*}
and
\begin{equation}
  \label{IsomorphismSobreGammaomega}
  \HH^{q}_{\dd'}(\cinfty_{\Gamma_{\pmb{\omega}}}(\Omega))
  \cong
  \HH^{q}_{\dd_t}(\cinfty_{\Gamma_{\pmb{\omega}}}(\Omega))
  \cong
  \cinfty(\TT^r) \otimes H^q_{\mathrm{dR}}(M);
\end{equation}
recall that
$r$ is the dimension of
$\Gamma_{\pmb{\omega}}$ as a
$\Z$-module.
We finish this section
with the following theorem, 
which introduces a condition
that will appear later in a
more general situation,
and that we will prove to hold
in several cases
(for notation,
see~\eqref{eq:cohomology-spaces-dprimexi}).
\begin{Thm}
  \label{thm:strong_isomorphism}
  Suppose that
  $\dd': \cinfty(\Omega; \LLambda^{q - 1}) \to \cinfty(\Omega; \LLambda^{q})$
  has closed range
  and
  \begin{equation}
    \label{eq:vanishing_mysterious_cohomologies}
    \HH_\xi^q (\cinfty(M)) = 0,
    \quad \forall \xi \in \Z^m \setminus \Gamma_{\pmb{\omega}}.
  \end{equation}
  Then
  \begin{equation*}
    \HH^{q}_{\dd'}(\cinfty_{\Z^m \setminus \Gamma_{\pmb{\omega}}}(\Omega)) = \{0\}.
  \end{equation*}
  In particular,
  in that case~\eqref{eq:dm_iso_intro} holds.
\end{Thm}
\begin{proof}
  Let
  $f \in \cinfty_{\Z^m \setminus \Gamma_{\pmb{\omega}}}(\Omega; \LLambda^{q})$
  be such that
  $\dd' f = 0$.
  Hence
  $\dd'_\xi \hat{f}_\xi = 0$,
  and by~\eqref{eq:vanishing_mysterious_cohomologies},
  there exists
  $u_\xi \in \cinfty(M; \Lambda^{q - 1})$
  such that
  $\dd'_\xi u_\xi = \hat{f}_\xi$
  for every
  $\xi \in \Z^m \setminus \Gamma_{\pmb{\omega}}$.
  Therefore
  $f$ is formally solvable,
  and by
  Theorem~\ref{Thm:reduction-glob-solv},
  there exists
  $u \in \cinfty(\Omega; \LLambda^{q - 1})$
  such that
  $\dd' u = f$.
  Replacing
  $u$ by
  $u_{\Z^m \setminus \Gamma_{\pmb{\omega}}}$
  yields the desired result.
\end{proof}

\subsection{Global solvability
  in the first degree}

Theorem~\ref{Thm:reduction-glob-solv}
tells us that
the obstruction to global solvability of
$\dd'$
is encoded in the frequencies
$\xi \in \Z^m \setminus \Gamma_{\pmb{\omega}}$.
This fact motivates us
to consider the following
Diophantine condition.
\begin{Def}
  \label{def:forms-simult-approx}
  A collection
  $\pmb{\omega} = (\omega_1, \ldots, \omega_m)$
  of real closed
  $1$-forms on
  $M$ is said to be
  \emph{strongly simultaneously approximable}
  if there exist a sequence
  of closed integral
  $1$-forms
  $\{ \beta_\nu \}_{\nu \in \N} \sset \cinfty(M; \Lambda^1)$
  and
  $\{ \xi_\nu \}_{\nu \in \N} \sset \Z^m \setminus \Gamma_{\pmb{\omega}}$
  such that
  $|\xi_\nu| \to \infty$ and
  \begin{equation}
    \label{eq:SSA_omega}
    \{ | \xi_\nu |^\nu (\xi_\nu \cdot \pmb{\omega} - \beta_\nu) \}_{\nu \in \N}
    \quad
    \text{is bounded in $\cinfty(M; \Lambda^1)$}.
  \end{equation}
  Otherwise,
  it is said to be
  \emph{weakly non-simultaneously approximable}.
\end{Def}
Such notions
depend only on the classes
$[\omega_1], \ldots, [\omega_m] \in H^1_{\mathrm{dR}} (M)$.
Indeed,
if for each
$k \in \{1, \ldots, m\}$
we have
$\omega_k^\bb \dfn \omega_k + \dd \gamma_k$
for some
$\gamma_k \in \cinfty(M; \R)$
then
\begin{equation*}
  \xi \cdot \pmb{\omega}^\bb
  =
  \xi \cdot \pmb{\omega}
  +
  \xi \cdot \dd \pmb{\gamma},
  \quad \forall \xi \in \Z^m,
\end{equation*}
in which
$\pmb{\omega}^\bb \dfn (\omega_1^\bb, \ldots, \omega_m^\bb)$
and
$\pmb{\gamma} \dfn (\gamma_1, \ldots, \gamma_m)$:
by integrating both sides against
an arbitrary
$1$-cycle,
it follows that
$\xi \cdot \pmb{\omega}^\bb$ is integral
if and only if so is
$\xi \cdot \pmb{\omega}$,
that is,
$\Gamma_{\pmb{\omega}^\bb} = \Gamma_{\pmb{\omega}}$.
Moreover,
\begin{equation*}
  | \xi_\nu |^\nu (\xi_\nu \cdot \pmb{\omega} - \beta_\nu)
  =
  | \xi_\nu |^\nu (\xi_\nu \cdot \pmb{\omega}^\bb - (\beta_\nu + \xi_\nu \cdot \dd \pmb{\gamma}))
\end{equation*}
and since each
$\beta_\nu + \xi_\nu \cdot \dd \pmb{\gamma}$
is integral,
it follows that 
$\pmb{\omega}^\bb$
is strongly simultaneously approximable
whenever
$\pmb{\omega}$ is.

As a finitely generated Abelian group,
$H_1(M; \Z)$ admits a primary decomposition
$\Z^d \oplus \Z_{p_1} \oplus \ldots \oplus \Z_{p_t}$.
We consider
$\sigma_1, \ldots, \sigma_d$ smooth
$1$-cycles whose homological classes
form a basis of its free part
$\Z^d$,
and the map
$I: H^1_{\mathrm{dR}}(M; \R) \to \R^d$
given by
\begin{equation*}
  I([\alpha])
  \dfn
  \frac{1}{2\pi} \left( \int_{\sigma_1} \alpha, \ldots, \int_{\sigma_d} \alpha \right).
\end{equation*}
Since
$\int_\sigma \alpha = 0$ if
$\sigma$ belongs to some
$\Z_{p_i}$,
it follows that a closed
$1$-form
$\alpha$ is integral
(resp.~rational)
if and only if
$I([\alpha]) \in \Z^d$
(resp.~$I[\alpha] \in \Q^d$).
By de Rham Theorem,
we conclude that
$I$ is a linear isomorphism
and that the classes of
$\sigma_1, \ldots, \sigma_d$,
regarded as elements of
$H_1(M; \R) \cong H^1_{\mathrm{dR}}(M; \R)^*$,
form an
$\R$-basis for the latter
vector space.

We associate to
$\pmb{\omega}$
the following family of vectors:
\begin{equation}
  \label{eq:v_ell}
  v_\ell
  \dfn
  \frac{1}{2 \pi}
  \left( \int_{\sigma_\ell} \omega_1, \ldots, \int_{\sigma_\ell} \omega_m \right) \in \R^m,
  \quad \ell \in \{1, \ldots, d \}.
\end{equation}
These, as well,
depend only on
$[\omega_1], \ldots, [\omega_m] \in H^1_{\mathrm{dR}} (M)$
by Stokes Theorem.
\begin{Prop}
  \label{Prop:non-simult-appr-charac}
  The collection
  $\pmb{\omega}$
  is strongly simultaneously approximable
  if and only if
  there are sequences
  $\{ \eta_\nu \}_{\nu \in \N} \sset \Z^d$
  and
  $\{ \xi_\nu \}_{\nu \in \N} \sset \Z^m \setminus \Gamma_{\pmb{\omega}}$
  such that
  $|\xi_\nu| \to \infty$
  and
  \begin{equation}
    \label{eq:SSA_vectors}
    \{ | \xi_\nu |^\nu (\xi_\nu \cdot v_\ell - \eta_{\nu \ell}) \}_{\nu \in \N}
    \quad
    \text{is bounded in $\R$
      for every
      $\ell \in \{ 1, \ldots, d \}$}.
  \end{equation}
\end{Prop}
\begin{proof}
  Let
  $Z^1 \sset \cinfty(M; \Lambda^1)$
  be the space of real, closed
  $1$-forms,
  endowed with
  the subspace topology,
  and define
  $J: Z^1 \to \R^d$
  by
  $J(\alpha) \dfn I([\alpha])$.
  It is continuous --
  as the composition of
  $I$ (a linear map between
  finite dimensional spaces)
  and the projection
  $Z^1 \to H^1_{\mathrm{dR}}(M; \R)$ --
  hence maps bounded sets to bounded sets.
  If
  $\pmb{\omega}$ is
  strongly simultaneously approximable,
  there exist a sequence
  of closed integral
  $1$-forms
  $\{ \beta_\nu \}_{\nu \in \N} \sset \cinfty(M; \Lambda^1)$
  and
  $\{ \xi_\nu \}_{\nu \in \N} \sset \Z^m \setminus \Gamma_{\pmb{\omega}}$
  such that
  $|\xi_\nu| \to \infty$
  such that~\eqref{eq:SSA_omega}
  holds.
  Letting
  $\eta_\nu \dfn I([\beta_\nu]) \in \Z^d$,
  we have,
  for each
  $\ell \in \{1, \ldots, d\}$,
  that
  \begin{equation}
    \label{eq:from_bdd_to_bdd}
    J ( | \xi_\nu |^\nu (\xi_\nu \cdot \pmb{\omega} - \beta_\nu) )_\ell
    =
    |\xi_\nu|^\nu (  J(\xi_\nu \cdot \pmb{\omega}) - J(\beta_\nu) )_\ell
    =
    |\xi_\nu|^\nu (\xi \cdot v_\ell - \eta_{\nu \ell})
  \end{equation}
  must be the
  $\ell$-th coordinate of
  a bounded sequence in
  $\R^d$,
  yielding our first claim.
    
  For the converse,
  we start with a digression.
  Let
  $\ker \Delta \sset Z^1$
  denote the space of harmonic
  $1$-forms w.r.t.~some
  Riemannian metric, which is well-known to be finite dimensional and
  hence has a well-defined norm topology;
  we have more, the map
  $\alpha \in \ker \Delta \mapsto [\alpha] \in H^1_{\mathrm{dR}}(M; \R)$
  is a linear isomorphism.
  Moreover,
  $\ker \Delta$ inherits a
  Fr{\'e}chet topology
  from
  $\cinfty(M; \Lambda^1)$,
  which matches the former one,
  by~\cite[Theorem~9.1]{treves_tvs}.
  Therefore the restriction
  $J: \ker \Delta \to \R^d$ 
  is a linear isomorphism,
  hence a topological one
  by finiteness;
  in particular, an
  $S \subset \ker \Delta$
  is bounded if and only if
  $J(S) \subset \R^d$ is bounded.
  We further pick
  $\vartheta_1, \ldots, \vartheta_d$
  a basis of
  $\ker \Delta$ dual to
  $[\sigma_1], \ldots, [\sigma_d] \in H_1(M; \R)$,
  in the sense that
  \begin{equation*}
    \frac{1}{2 \pi} \int_{\sigma_\ell} \vartheta_{\ell'}
    =
    \delta_{\ell \ell'},
    \quad \forall \ell, \ell' \in \{1, \ldots, d\}.
  \end{equation*}
  Such a choice makes
  $\vartheta_1, \ldots, \vartheta_d$
  integral.
  
  Now, suppose there are sequences
  $\{ \eta_\nu \}_{\nu \in \N} \sset \Z^d$
  and
  $\{ \xi_\nu \}_{\nu \in \N} \sset \Z^m \setminus \Gamma_{\pmb{\omega}}$
  such that
  $|\xi_\nu| \to \infty$
  and~\eqref{eq:SSA_vectors}
  holds.
  Then
  \begin{equation*}
    \beta_\nu
    \dfn
    \sum_{\ell = 1}^d \eta_{\nu \ell} \vartheta_\ell
    \in \ker \Delta
  \end{equation*}
  is integral and
  $J(\beta_\nu) = \eta_\nu$
  for every
  $\nu \in \N$.
  Since~\eqref{eq:from_bdd_to_bdd}
  holds once more
  for each
  $\ell \in \{1, \ldots, d\}$,
  we have from~\eqref{eq:SSA_vectors}
  that
  $\{ J (| \xi_\nu |^\nu (\xi_\nu \cdot \pmb{\omega} - \beta_\nu) ) \}_{\nu \in \N}$
  must be a bounded sequence in
  $\R^d$:
  the conclusion follows from
  our digression,
  since we may assume
  w.l.o.g.~$\omega_1, \ldots, \omega_m \in \ker \Delta$.
\end{proof}

\begin{Cor}
  \label{cor:non-simult-appr-charac}
  The collection
  $\pmb{\omega}$
  is weakly non-simultaneously approximable
  if and only if there exist
  $C, \rho > 0$
  such that
  (see~\eqref{eq:v_ell})
  \begin{equation}
    \label{eq:simult-appr}
    \max_{1 \leq \ell \leq d} | \xi \cdot v_\ell - \eta_\ell |
    \geq
    C |\xi|^{-\rho},
    \quad \forall \eta \in \Z^d, \ \forall \xi \in \Z^m \setminus \Gamma_{\pmb{\omega}}.
  \end{equation}
\end{Cor}
\begin{proof}
  Suppose~\eqref{eq:simult-appr}
  does not hold for any
  $C,\rho > 0$.
  Then, for every
  $\nu \in \N$, we can find
  $\eta_\nu \in \Z^d$ and
  $\xi_\nu \in \Z^m \setminus \Gamma_{\pmb{\omega}}$
  such that
  \begin{equation}
    \label{eq:aux_SSA1}
    |\xi_\nu|^\nu | \xi_\nu \cdot v_\ell - \eta_{\nu \ell} | < \nu^{-1},
    \quad \forall \nu \in \N, \ \forall \ell \in \{1, \ldots, d\}.
  \end{equation}
  We claim that the sequence
  $\{ \xi_\nu \}_{\nu \in \N}$ is unbounded,
  otherwise, we would have
  \begin{equation*}
    |\eta_{\nu \ell}|
    \leq
    | \xi_\nu \cdot v_\ell - \eta_{\nu \ell} | + | \xi_\nu| |v_\ell|
    \Longrightarrow
    |\eta_\nu| \leq \nu^{-1} |\xi_\nu|^{-\nu} + | \xi_\nu| \max_{1 \leq \ell \leq d} |v_\ell|,
  \end{equation*}
  hence
  $\{ \eta_\nu \}_{\nu \in \N}$
  would be also bounded,
  in which case
  $\{ (\eta_\nu, \xi_\nu) \}_{\nu \in \N}$
  would attain at most finitely many distinct values,
  making~\eqref{eq:aux_SSA1} impossible.
  This proves our claim,
  and by extracting an increasing subsequence of
  $\{ \xi_\nu \}_{\nu \in \N}$,
  we conclude that
  $\pmb{\omega}$
  is strongly simultaneously approximable
  by means of
  Proposition~\ref{Prop:non-simult-appr-charac}.

  Conversely,
  suppose that~\eqref{eq:simult-appr}
  holds for some
  $C,\rho > 0$,
  and, by contradiction, that
  $\pmb{\omega}$
  is strongly simultaneously approximable.
  By Proposition~\ref{Prop:non-simult-appr-charac},
  there exist
  $C' > 0$, and sequences
  $\{ \eta_\nu \}_{\nu \in \N} \sset \Z^d$
  and
  $\{ \xi_\nu \}_{\nu \in \N} \sset \Z^m \setminus \Gamma_{\pmb{\omega}}$
  with
  $|\xi_\nu| \to \infty$,
  such that
  \begin{equation*}
    C' |\xi_\nu|^{-\nu}
    \geq
    \max_{1 \leq \ell \leq d} | \xi_\nu \cdot v_\ell - \eta_{\nu \ell} |
    \geq
    C |\xi_\nu|^{-\rho}
    \Longrightarrow
    C^{-1} C' \geq |\xi_\nu|^{\nu - \rho} \lra \infty,
  \end{equation*}
  an evident contradiction.
\end{proof}
\begin{Rem}
  \hfill
  \begin{enumerate}
  \item Notice that
    \begin{equation*}
      \Gamma_{\pmb{\omega}}
      =
      \{ \xi \in \Z^m \st \xi \cdot v_\ell \in \Z,
      \ \forall \ell \in \{1, \ldots, d\} \};
    \end{equation*}
  \item The inequality~\eqref{eq:simult-appr}
    is equivalent to
    condition~$\mathrm{(DC)}$
    in~\cite[Section~2]{dm16}
    for the
    $d \times m$ matrix whose
    rows are
    $v_1, \ldots, v_d$.
    When 
    $\Gamma_{\pmb{\omega}} = \{0\}$,
    it recovers the standard notion
    of non-simultaneous approximability
    for collections of vectors
    in~\cite[Definition~1.1]{hp00}
    (see further connections with previous conditions
    in the literature there);
  \item If
    $v_1, \ldots, v_d \in \Q^m$,
    which corresponds to case when
    $\omega_1, \ldots, \omega_m$
    are all rational
    $1$-forms,
    we pick a non-zero
    $\lambda \in \Z_+$
    such that
    $\lambda v_\ell \in \Z^m$
    for every
    $\ell \in \{1, \ldots, d\}$.
    Hence, given
    $\xi \in \Z^m \setminus \Gamma_{\pmb{\omega}}$,
    there must exist an
    $\ell_0 \in \{1, \ldots, d\}$
    for which
    $\xi \cdot v_{\ell_0} \notin \Z$,
    thus, for any
    $\eta \in \Z^d$, we must have
    \begin{align*}
      \xi \cdot v_{\ell_0} - \eta_{\ell_0} \neq 0
      &\Longrightarrow
      \lambda (\xi \cdot v_{\ell_0} - \eta_{\ell_0}) \in \Z \setminus 0 \\
      &\Longrightarrow
      \lambda \max_{1 \leq \ell \leq d} | \xi \cdot v_\ell - \eta_\ell | \geq 1,
    \end{align*}
    that is,~\eqref{eq:simult-appr} holds with
    $C \dfn \lambda^{-1}$
    and
    $\rho > 0$ arbitrary.
  \end{enumerate}
\end{Rem}

Now, we can state
our characterization of solvability
in the first degree
in terms of the
Diophantine condition
just introduced.
\begin{Thm}
  \label{Thm:solv-diop}
  The map
  $\dd': \cinfty(\Omega) \to \cinfty(\Omega; \LLambda^{1})$
  has closed range
  if and only if the collection
  $\pmb{\omega}$
  is weakly non-simultaneously approximable.
\end{Thm}
\begin{proof}
  $(\Longrightarrow)$: Suppose  there exist
  a sequence of integral
  $1$-forms
  $\{ \beta_\nu \}_{\nu \in \N} \sset \cinfty(M; \Lambda^1)$
  and
  $\{ \xi_\nu \}_{\nu \in \N} \sset \Z^m \setminus \Gamma_{\pmb{\omega}}$
  such that
  $|\xi_\nu| \to \infty$
  and~\eqref{eq:SSA_omega}
  holds.
  We may assume that
  $\xi_\nu \neq \xi_{\nu'}$
  whenever
  $\nu \neq \nu'$.
  Let
  $\varphi_\nu \in \cinfty (\tilde M; \R)$
  be such that
  $\dd \varphi_\nu = \pi^* \beta_\nu$
  and define
  \begin{equation*}
    f(t,x)
    \dfn
    \frac{1}{(2 \pi)^m} \sum_{\nu = 1}^\infty
    e^{ix\xi_\nu-i \varphi_\nu \circ \pi^{-1}(t)}  (\xi_\nu \cdot \pmb{\omega} - \beta_\nu)(t),
  \end{equation*}
  which, we claim,
  is smooth and formally solvable.
  As we saw
  in~\eqref{estabemdefinidaemM},
  we have that
  \begin{equation*}
    M \ni t \longmapsto e^{-i \varphi_\nu \circ \pi^{-1}(t)}
  \end{equation*}
  is a well-defined smooth function since
  $\beta_\nu$ is integral.
  Therefore, we define, for every
  $\xi \in \Z^m$,
  a smooth $1$-form on $M$
  \begin{equation*}
    f_\xi
    \dfn
    \begin{cases}
      e^{-i \varphi_\nu \circ \pi^{-1}} (\xi_\nu \cdot \pmb{\omega} - \beta_\nu),
      &\text{if $\xi = \xi_\nu$}, \\
      0,
      &\text{otherwise}.
    \end{cases}
  \end{equation*}
  Note that, in order to prove that
  $f$ is a well-defined smooth
  $1$-form, it is enough to prove
  that every point of
  $M$ belongs to a coordinate system
  $(U; t_1, \ldots, t_n)$,
  with
  $\alpha \in \Z_+^n$
  and
  $s > 0$,
  there exists
  $C > 0$
  such that
  \begin{equation*}
    \sup_U \| \partial_t^\alpha f_{\xi_\nu} \|
    \leq
    C (1 + |\xi_\nu|)^{-s},
    \quad \forall \nu \in \N,
  \end{equation*}
  in which the norm
  $\| \cdot \|$ in the left-hand side is
  the sum of the absolute values
  of the coordinate components of the
  $1$-form
  w.r.t.~the local frame
  $\dd t_1, \ldots, \dd t_n$.
  By hypothesis~\eqref{eq:SSA_omega},
  we may assume that, for each
  $\gamma \in \Z_+^n$,
  there exists
  $B > 0$
  such that
  \begin{equation*}
    \sup_M \| \partial_t^\gamma (\xi_\nu \cdot \pmb{\omega} - \beta_\nu) \|
    \leq
    B |\xi_\nu|^{-\nu},
    \quad \forall \nu \in \N,
  \end{equation*}
  It remains to prove that
  the derivatives of
  $e^{-i\varphi_\nu \circ \pi^{-1}}$
  are bounded by some constant times a power of
  $|\xi_\nu|$.
  We can assume that
  $U$ is small enough in order to
  have a diffeomorphism
  $\pi_U: \tilde{U} \to U$
  and so
  \begin{equation*}
    e^{-i\varphi_\nu \circ \pi^{-1}(t)} = e^{-i\varphi_\nu \circ \pi^{-1}_U(t)}  
  \end{equation*}
  for every
  $t \in U$ and so,
  since
  $\varphi_\nu$ is real,
  our claim follows from the identity
  \begin{equation*}
    \sum_{j = 1}^n \frac{\partial (\varphi_\nu \circ \pi_U^{-1})}{\partial t_j} \dd t_j
    =
    \beta_\nu
    =
    - (\xi_\nu \cdot \pmb{\omega} - \beta_\nu) + \xi_\nu \cdot \pmb{\omega},
  \end{equation*}
  since the right-hand
  of the equality above
  clearly is bounded by a constant times
  $|\xi_\nu|$,
  hence proving that
  $f \in \cinfty(\Omega; \LLambda^{1})$.
  The next step
  is to prove that
  $f$ is formally solvable.
  Since
  $\hat f_\xi = 0$, if
  $\xi \notin \{\xi_\nu \}_{\nu \in \N}$,
  it follows that
  $\hat f_\xi$ is
  $\dd'_\xi$-exact
  for such values of
  $\xi$;
  for
  $\xi = \xi_\nu$
  we have
  \begin{equation*}
    \dd'_{\xi_\nu} e^{-i \varphi_\nu\circ \pi^{-1}}
    =
    \dd e^{-i \varphi_\nu\circ \pi^{-1}} + i e^{-i \varphi_\nu\circ \pi^{-1}} (\xi_\nu \cdot \pmb{\omega})
    =
    -i e^{-i \varphi_\nu\circ \pi^{-1}} \beta_\nu + i e^{-i \varphi_\nu\circ \pi^{-1}} (\xi_\nu \cdot \pmb{\omega})
    =
    \hat{f}_{\xi_\nu}.
  \end{equation*}
  
  It remains to prove that
  there is no
  $u \in \cinfty(\Omega)$
  such that
  $\dd' u = f$.
  If such an
  $u$ exists, then
  $\dd'_{\xi} \hat{u}_\xi = \hat f_{\xi}$
  for every
  $\xi \in \Z^m$,
  in particular
  $\hat u_{\xi_\nu} - e^{-i\varphi_\nu\circ \pi^{-1}} \in \ker \dd'_{\xi_\nu}$
  for every
  $\nu \in \N$.
  But since
  $\xi_\nu \notin \Gamma_{\pmb{\omega}}$,
  we have by~\cite[Lemma~2.1]{bcp96}
  that $\ker \dd'_{\xi_\nu} = \{0\}$,
  hence
  $\hat u_{\xi_\nu} = e^{-i\varphi_\nu\circ \pi^{-1}}$
  for every
  $\nu \in \N$
  and then
  \begin{equation*}
    \| \hat u_{\xi_\nu} \|_{L^2(M)}^2  = \int_M \dd V_M,
    \quad \nu \in \N,
  \end{equation*}
  does not decrease fast,
  contradicting the smoothness of
  $u$.

  $(\Longleftarrow)$: Suppose that
  $f \in \cinfty(\Omega ;\LLambda^{1})$
  is formally solvable,
  i.e., for every
  $\xi \in \Z^m$, there exists
  $u_\xi \in \cinfty(M)$
  such that
  $\dd'_\xi u_\xi = \hat f_\xi$.
  Thanks to
  Theorem~\ref{Thm:reduction-glob-solv},
  we only need to prove that
  there exists
  $u \in \cinfty(\Omega)$
  such that
  $\dd' u = f_{\Z^m\setminus \Gamma_{\pmb{\omega}}}$.
  The natural candidate
  is 
  \begin{equation*}
    u
    \dfn
    \frac{1}{(2\pi)^m} \sum_{\xi \in \Z^m \setminus \Gamma_{\pmb{\omega}}} e^{ix \xi} u_\xi,
  \end{equation*}
  and it suffices to prove that
  the series above converges in
  $\cinfty(\Omega)$.
  This claim is local,
  so we must verify that
  every point of
  $M$ belongs to a coordinate ball
  $U$
  enjoying the following property:
  for each
  $\alpha \in \Z_+^n$
  and
  $s > 0$,
  there exists
  $C > 0$
  such that
  \begin{equation}
    \label{Eq:est-fourier-part-sol}
    | \partial_t^\alpha u_\xi(t) |
    \leq C(1 + |\xi|)^{-s},
    \quad \forall \xi \in \Z^m \setminus \Gamma_{\pmb{\omega}}, \ \forall t \in U.
  \end{equation}
  We write
  \begin{equation*}
    \omega_k = \sum_{j = 1}^n \omega_{kj} \dd t_j,
    \quad \omega_{kj} \in \cinfty(U),
  \end{equation*}
  on $U$
  for
  $k \in \{1,\ldots,m\}$,
  and
  \begin{equation*}
    f = \sum_{j = 1}^n f_j \ \dd t_j,
    \quad f_j \in \cinfty(U \times \TT^m),
  \end{equation*}
  on
  $U \times \TT^m$.
  We may assume
  $\omega_{kj}$ is bounded in the
  $\sup$ norm in $U$.
  \begin{Lem}
    \label{Lem:est-fourier-solution}
    Suppose that given
    $s > 0$
    there exists
    $C > 0$
    such that
    \begin{equation}
      \label{Eq:est-fourier-part-sol'}
      | u_\xi(t) |
      \leq
      C(1 + |\xi|)^{-s},
      \quad \forall \xi \in \Z^m\setminus \Gamma_{\pmb{\omega}}, \quad  \forall t \in U.
    \end{equation}
    Then~\eqref{Eq:est-fourier-part-sol}
    holds true.
  \end{Lem} 
  \begin{proof}
    We proceed by induction on
    $|\alpha|$;
    the base
    $|\alpha| = 0$
    is~\eqref{Eq:est-fourier-part-sol'}.
    For the general case,
    the equality
    $\dd'_\xi u_\xi = \hat f_\xi$
    can be written in
    $U$ as
    \begin{equation*}
      \frac{\partial u_\xi}{\partial t_j}
      + i \sum_{k = 1}^m \xi_k \omega_{kj} u_\xi
      =
      \mathcal{F}_\xi(f_j),
      \quad j \in \{1, \ldots, n\}.
    \end{equation*}
    If
    $\alpha = (\alpha_1, \ldots, \alpha_n)$,
    $\alpha_j > 0$
    and
    $\beta \dfn \alpha - e_j$,
    in which
    $e_j$ is the
    $j$-th vector of the canonical basis of
    $\R^n$,
    then
    \begin{equation*}
      \partial_t^\alpha u_\xi
      =
      \partial_t^\beta \partial_{t_j} u_\xi
      = \partial_t^\beta
      \left( \mathcal{F}_\xi (f_j) - i \sum_{k = 1}^m \xi_k \omega_{kj} u_\xi \right),
    \end{equation*}
    and we use that
    $f$ is smooth and
    $s > 0$
    in~\eqref{Eq:est-fourier-part-sol}
    is arbitrary.
  \end{proof}

  We fix
  $t_0 \in M$
  and prove
  that~\eqref{Eq:est-fourier-part-sol'}
  holds in a coordinate ball
  $U$ centered at
  $t_0$.
  As in Section~\ref{sec:isomorphisms},
  we assume that
  $\pi_U: \tilde U \to U$
  is a diffeomorphism
  and consider the functions
  $\psi_\xi \in \cinfty(\tilde{M}; \R)$.
  Then on
  $U$ we have,
  for all
  $\xi \in \Z^m$,
  \begin{equation*}
    \dd (e^{i \psi_\xi \circ \pi_U^{-1}} u_\xi )
    =
    e^{i \psi_\xi \circ \pi_U^{-1}} \dd'_\xi u_\xi
    =
    e^{i \psi_\xi \circ \pi_U^{-1}} \hat{f}_\xi
  \end{equation*}
  which,
  integrating along any curve
  $\gamma$ in $U$ connecting
  $t_0$ to $t$ yields,
  by Stokes Theorem,
  \begin{equation*}
    u_\xi(t)
    =
    e^{i(\psi_\xi \circ \pi_U^{-1}(t_0) - \psi_\xi \circ \pi_U^{-1}(t))} u_\xi(t_0)
    +
    e^{-i\psi_\xi \circ \pi_U^{-1}(t)} \int_{\gamma} e^{i\psi_\xi\circ \pi_U^{-1}} \hat f_\xi.
  \end{equation*}
  Since
  $f$ is smooth and 
  $\psi_\xi$ is real-valued,
  the last equality shows
  that~\eqref{Eq:est-fourier-part-sol'}
  holds true provided that we prove that, given
  $s > 0$, there exists
  $B > 0$ such that
  \begin{equation*}
    | u_\xi(t_0) | \leq B (1 + |\xi|)^{-s},
    \quad \forall \xi \in \Z^m\setminus \Gamma_{\pmb{\omega}},
  \end{equation*}
  which we do next.

  By the Hurewicz's Theorem,
  we may assume that 
  $\sigma_1, \ldots, \sigma_d$
  all have base point
  $t_0$.
  We lift them to curves
  $\tilde \sigma_\ell: [0,1] \to \tilde M$
  whose endpoints we denote by
  \begin{equation*}
    Q_0 \dfn \tilde \sigma_\ell(0),
    \quad
    Q_\ell \dfn \tilde \sigma_\ell(1),
    \quad \ell \in \{1, \ldots, d\}.
  \end{equation*}
  We integrate along
  $\tilde \sigma_\ell$
  both sides of the equality
  \begin{equation*}
    \dd (e^{i \psi_\xi} \pi^* u_\xi )
    =
    e^{i \psi_\xi} \pi^* \hat f_\xi
  \end{equation*}
  and use Stokes Theorem again
  to conclude that,
  since
  $\pi(Q_\ell) = \pi(Q_0) = t_0$,
  \begin{equation*}
    \left( e^{i (\psi_\xi(Q_\ell) - \psi_\xi(Q_0))} - 1\right) u_\xi(t_0)
    =
    e^{i\psi_\xi(Q_0)} \int_{\tilde \sigma_\ell} e^{i\psi_\xi} \pi^* \hat f_\xi.
  \end{equation*}
  Moreover,
  \begin{equation}
    \label{eq:formula-integral-psi}
    \psi_\xi(Q_\ell) - \psi_\xi(Q_0)
    =
    \int_{\partial \tilde \sigma_\ell} \psi_\xi
    =
    \int_{\tilde \sigma_\ell} \dd \psi_\xi
    =
    \int_{\tilde \sigma_\ell} \pi^* (\xi \cdot \pmb{\omega})
    =
    \int_{\sigma_\ell} \xi \cdot \pmb{\omega},
  \end{equation}
  so, for
  $\xi \in \Z^m \setminus \Gamma_{\pmb{\omega}}$,
  there exists
  $\ell$
  such that
  $\psi_\xi(Q_\ell)-\psi_\xi(Q_0) \notin 2 \pi \Z$,
  hence
  \begin{equation*}
    u_\xi(t_0)
    =
    \frac{e^{i\psi_\xi(Q_0)}}
         {\left( e^{i (\psi_\xi(Q_\ell) - \psi_\xi(Q_0))} - 1 \right)}
         \int_{\tilde \sigma_\ell} e^{i\psi_\xi} \pi^* \hat f_\xi
  \end{equation*}
  for any such
  $\ell$.
  Using again that
  $f$ is smooth and
  $\psi_\xi$ is real-valued,
  it suffices to prove that
  there exist
  $c, \rho > 0$
  such that
  \begin{equation}
    \label{eq:est-exp}
    \max_{1 \leq \ell \leq d}
    \left| e^{i (\psi_\xi(Q_\ell) - \psi_\xi(Q_0))} - 1 \right|
    \geq
    c |\xi|^{-\rho},
    \quad \forall \xi \in \Z^m \setminus \Gamma_{\pmb{\omega}}.
  \end{equation}
  Let us assume for a moment the following lemma:
  \begin{Lem}
    \label{Lem:ineq-aux-exp}
    There exists
    $\epsilon > 0$
    such that, for each
    $\xi \in \Z^m \setminus \Gamma_{\pmb{\omega}}$,
    at least one of the
    following conditions holds:
    \begin{enumerate}
    \item For every
      $\ell \in \{1, \ldots, d\}$
      there exists
      $p_\ell \in \Z$
      such that
      \begin{equation*}
        \left|e^{i(\psi_\xi(Q_0) - \psi_\xi(Q_\ell))} - 1 \right|
        \geq
        \frac{1}{2}
        \left| \psi_\xi(Q_0) - \psi_\xi(Q_\ell) - 2 \pi p_\ell \right|;
      \end{equation*}
    \item There exists
      $\ell \in \{1, \ldots, d\}$
      such that
      \begin{equation*}
        \left| e^{i (\psi_\xi(Q_0) - \psi_\xi(Q_\ell) )} - 1 \right|
        \geq \epsilon.
      \end{equation*}
    \end{enumerate}
  \end{Lem}
  By Lemma~\ref{Lem:ineq-aux-exp},
  we may assume that,
  for every
  $\xi \in \Z^m \setminus \Gamma_{\pmb{\omega}}$,
  there are
  $p_1, \ldots, p_d \in \Z$
  such that
  \begin{equation*}
    \max_{1 \leq \ell \leq d}
    \left| e^{i (\psi_\xi(Q_\ell) - \psi_\xi(Q_0))} - 1 \right|
    \geq
    \frac{1}{2} \max_{1 \leq \ell \leq d}
    \left| \psi_\xi(Q_\ell) - \psi_\xi(Q_0) - 2\pi p_\ell  \right|.
  \end{equation*}
  It follows
  from~\eqref{eq:formula-integral-psi}
  and~\eqref{eq:v_ell}
  that
  \begin{equation*}
    \max_{1 \leq \ell \leq d}
    \left| e^{i (\psi_\xi(Q_\ell) - \psi_\xi(Q_0))} - 1 \right|
    \geq
    \pi \max_{1 \leq \ell \leq d} | \xi \cdot v_\ell - p_\ell |.
  \end{equation*}
  Since we are assuming that
  $\pmb{\omega}$ is weakly non-simultaneously approximable,
  we can use
  Corollary~\ref{cor:non-simult-appr-charac}
  to conclude
  that~\eqref{eq:est-exp} holds.
\end{proof}
\begin{proof}[Proof of Lemma~\ref{Lem:ineq-aux-exp}]
  Given
  $\xi \in \Z^m \setminus \Gamma_{\pmb{\omega}}$,
  exactly one alternative holds:
  \begin{enumerate}
  \item For each
    $\ell \in \{1, \ldots, d\}$
    there exists
    $p_\ell \in \Z$ such that
    \begin{equation}
      \label{eq:ineq-aux-exp1}
      \left| \psi_\xi(Q_0) - \psi_\xi(Q_\ell) - 2 \pi p_\ell \right|
      <
      \frac{1}{2};
    \end{equation}
  \item
    There exists
    $\ell \in \{1, \ldots, d\}$
    such that
    \begin{equation}
      \label{eq:ineq-aux-exp2}
      \left| \psi_\xi(Q_0) - \psi_\xi(Q_\ell) - 2 \pi p \right|
      \geq
      \frac{1}{2},
      \quad \forall p \in \Z.
    \end{equation}
  \end{enumerate}
  Take
  $\epsilon \dfn 1 - \cos (1/2)$.
  If~\eqref{eq:ineq-aux-exp2}
  holds for some
  $\ell$, then
  \begin{equation*}
    \left| e^{i(\psi_\xi(Q_0) - \psi_\xi(Q_\ell))} - 1 \right|
    \geq
    \left| \cos(\psi_\xi(Q_0) - \psi_\xi(Q_\ell)) - 1 \right|
    \geq
    \epsilon
  \end{equation*}
  since
  $\epsilon \leq 1 - \cos x$ for any
  $x \in [2p \pi +\frac{1}{2}, 2(p+1)\pi - \frac{1}{2}]$,
  for every
  $p \in \Z$.      
  
  Otherwise,
  if~\eqref{eq:ineq-aux-exp1}
  happens to be true,
  then, for
  $z \dfn i (\psi_\xi(Q_0) - \psi_\xi(Q_\ell))$,
  we have
  \begin{equation*}
    |e^{z} - 1 - (z - i2 \pi p_\ell)|
    =
    \left|
    \sum_{n = 2}^\infty
    \frac{(z - i2 \pi p_\ell)^{n}}{n!}
    \right|
    \leq
    |z - i2 \pi p_\ell|^2
    \sum_{n = 2}^\infty \frac{(1/2)^{n - 2}}{n!}
    \leq
    |z - i2 \pi p_\ell|^2.
  \end{equation*}
  Since
  $|a - b| \leq |b|^2$
  and
  $|b| \leq 1/2$
  imply
  $|b| /2 \leq |a|$,
  our conclusion
  follows at once.
\end{proof}

\section{Formal Fourier analysis of cohomology spaces}
\label{sec:fourier_cohom_spaces}

Allowing again
$\omega_1, \ldots, \omega_m$
to be complex-valued,
we add a third characterization
of solvability to our list in
Corollary~\ref{cor:formal_iff_closed}.
\begin{Prop}
  \label{prop:formal_iff_into_prod}
  The following are equivalent:
  \begin{enumerate}
  \item For every formally solvable
    $f \in \cinfty(\Omega; \LLambda^{q})$ 
    there exists
    $u \in \cinfty(\Omega; \LLambda^{q - 1})$
    such that
    $\dd' u = f$;
  \item The map
    \begin{equation}
      \label{eq:into_product}
      \HH^{q}_{\dd'}(\cinfty(\Omega)) \lra \prod_{\xi \in \Z^m} \HH^q_{\xi}(\cinfty(M))
    \end{equation}
    given by
    $[f] \mapsto ([\hat{f}_\xi])_{\xi \in \Z^m}$
    is injective.
  \end{enumerate}
\end{Prop}
\begin{proof}
  First, we prove that
  solvability implies the injectivity
  of~\eqref{eq:into_product}. 
  Let
  $[f] \in \HH^{q}_{\dd'}(\cinfty(\Omega))$
  be such that
  $[\hat{f}_\xi] = 0$ in
  $\HH^q_{\xi}(\cinfty(M))$
  for every $\xi \in \Z^m$,
  i.e., $\hat{f}_\xi$ is
  $\dd'_\xi$-exact
  for each $\xi \in \Z^m$,
  i.e., $f$ is formally solvable.
  By the solvability, it follows that
  $[f] = 0$ in
  $\HH^{q}_{\dd'}(\cinfty(\Omega))$,
  proving the injectivity. 

  Conversely, a formally solvable
  $f \in \cinfty(\Omega; \LLambda^{q})$
  is $\dd'$-closed,
  hence it determines a class
  $[f] \in \HH^{q}_{\dd'}(\cinfty(\Omega))$.
  Formal solvability of
  $f$ ensures that
  $[\hat{f}_\xi] = 0$ in
  $\HH^q_{\xi}(\cinfty(M))$
  for every $\xi \in \Z^m$. Thus, 
  if~\eqref{eq:into_product}
  is  injective,
  we have
  $[f] = 0$ in
  $\HH^{q}_{\dd'}(\cinfty(\Omega))$,
  which is precisely what we wanted to prove.
\end{proof}

Next, we study the map induced by
$\mathcal{E}_\xi$ in cohomology.
Unlike in the previous proposition
for the map induced by
$\mathcal{F}_\xi$,
its injectivity always holds true
and requires no extra hypotheses.
\begin{Prop}
  \label{prop:into_sum}
  For each
  $\xi \in \Z^m$,
  $\mathcal{E}_\xi$
  induces an injection
  $\HH_\xi^q (\cinfty(M)) \hookrightarrow \HH_{\dd'}^{q}(\cinfty(\Omega))$.
  Their direct sum induces
  the following injection:
  \begin{equation}
    \label{eq:into_sum}
    \bigoplus_{\xi \in \Z^m} \HH_\xi^q (\cinfty(M)) \lra \HH_{\dd'}^{q}(\cinfty(\Omega)).
  \end{equation}
\end{Prop}
\begin{proof}
  Let
  $[f] \in \HH_\xi^q (\cinfty(M))$.
  Then,
  by Corollary~\ref{cor:actions_fouriercoeff},
  \begin{equation*}
    \dd' \mathcal{E}_\xi f
    = \mathcal{E}_\xi \dd'_\xi f
    = 0
  \end{equation*}
  and the class of
  $\mathcal{E}_\xi f$ in
  $\HH_{\dd'}^{q}(\cinfty(\Omega))$
  is well-defined. Indeed,
  if
  $[f^\bb] = [f]$ in
  $\HH_\xi^q (\cinfty(M))$,
  then there exists
  $v \in \cinfty(M; \Lambda^{q - 1})$
  such that
  $\dd'_\xi v = f - f^\bb$,
  hence
  \begin{equation*}
    \mathcal{E}_\xi f - \mathcal{E}_\xi f^\bb
    = \mathcal{E}_\xi \dd'_\xi v
    = \dd' \mathcal{E}_\xi v.
  \end{equation*}
  Thus
  $[\mathcal{E}_\xi f^\bb] = [\mathcal{E}_\xi f]$
  in $\HH_{\dd'}^{q}(\cinfty(\Omega))$.
  Moreover, if
  $[\mathcal{E}_\xi f] = 0$ in
  $\HH_{\dd'}^{q}(\cinfty(\Omega))$,
  then there exists
  $u \in \cinfty(\Omega; \LLambda^{q - 1})$
  such that
  $\dd' u = \mathcal{E}_\xi f$,
  hence, by Lemma~\ref{lem:rightinverseFxiExi}, we have
  \begin{equation*}
    f
    = \mathcal{F}_\xi \mathcal{E}_\xi f
    = \mathcal{F}_\xi \dd' u
    = \dd'_\xi \hat{u}_\xi
  \end{equation*}
  hence the assignment
  $[f] \mapsto [\mathcal{E}_\xi f]$
  is injective.
  More generally,
  if
  \begin{equation*}
    \sum_{\xi \in \Z^m} [f_\xi] \in \bigoplus_{\xi \in \Z^m} \HH_\xi^q (\cinfty(M))
  \end{equation*}
  is such that
  \begin{equation*}
    \sum_{\xi \in \Z^m} \mathcal{E}_\xi f_\xi \in \cinfty(\Omega; \LLambda^{q})
    \quad
    \text{is $\dd'$-exact}
  \end{equation*}
  (the smoothness is possible by taking
  $0$ as representative whenever
  $[f_\xi] = 0$),
  then for some
  $u \in \cinfty(\Omega; \LLambda^{q - 1})$,
  we have
  \begin{equation*}
    \dd' u
    = \sum_{\xi \in \Z^m} \mathcal{E}_\xi f_\xi
    = \frac{1}{(2 \pi)^m} \sum_{\xi \in \Z^m} e^{ix \xi} \wedge f_\xi
    \Longrightarrow
    \dd'_\xi \hat{u}_\xi = f_\xi,
    \quad \forall \xi \in \Z^m,
  \end{equation*}
  that is,
  $[f_\xi] = 0$ in
  $\HH_\xi^q (\cinfty(M))$
  for every $\xi \in \Z^m$.
\end{proof}

\section{General finiteness theorems}
\label{sec:gen_fin_thm}

The next result
is the heart of our
forthcoming analysis.
\begin{Thm}
  \label{thm:finitedimensional}
  Given
  $q \in \{0, \ldots, n\}$,
  the following are equivalent:
  \begin{enumerate}
  \item $\HH_{\dd'}^{q}(\cinfty(\Omega))$
    is finite dimensional;
  \item $\dd': \cinfty(\Omega; \LLambda^{q - 1}) \to \cinfty(\Omega; \LLambda^{q})$
    has closed range and
    there exists a finite set
    $F \sset \Z^m$
    such that
    \begin{equation*}
      \HH_\xi^q (\cinfty(M)) = \{0\},
      \quad \forall \xi \in \Z^m \setminus F.
    \end{equation*}
  \end{enumerate}
  In this case:
  \begin{equation*}
    \HH_{\dd'}^{q}(\cinfty(\Omega)) \cong \bigoplus_{\xi \in F} \HH_\xi^q (\cinfty(M)).
  \end{equation*}
\end{Thm}
\begin{proof}
  It is a standard argument in
  Functional Analysis that, if
  $\HH_{\dd'}^{q}(\cinfty(\Omega))$
  is finite dimensional,
  then the denominator
  in~\eqref{eq:cohomology_VV}
  is a closed subspace of
  $\cinfty(\Omega; \LLambda^{q})$.
  Moreover, since
  the map~\eqref{eq:into_sum}
  is injective,
  finiteness of
  $\HH_{\dd'}^{q}(\cinfty(\Omega))$
  entails that only
  finitely many terms in
  that direct sum
  are non-zero.

  Conversely, if
  $\dd': \cinfty(\Omega; \LLambda^{q - 1}) \to \cinfty(\Omega; \LLambda^{q})$
  has closed range,
  then
  the map~\eqref{eq:into_product}
  is injective
  (by
  Corollary~\ref{cor:formal_iff_closed}
  and
  Proposition~\ref{prop:formal_iff_into_prod}),
  i.e.~$\HH_{\dd'}^{q}(\cinfty(\Omega))$
  injects into
  \begin{equation*}
    \prod_{\xi \in F} \HH^q_{\xi}(\cinfty(M))
  \end{equation*}
  which is
  finite dimensional
  since so is every factor and
  $F$ is finite.
  It follows that
  $\dim \HH_{\dd'}^{q}(\cinfty(\Omega)) < \infty$,
  and we have a sequence of
  injections~\eqref{eq:into_sum}-\eqref{eq:into_product}
  between finite dimensional spaces
  \begin{equation*}
    \bigoplus_{\xi \in F} \HH^q_{\xi}(\cinfty(M))
    \lra
    \HH_{\dd'}^{q}(\cinfty(\Omega))
    \lra
    \prod_{\xi \in F} \HH^q_{\xi}(\cinfty(M))
  \end{equation*}
  with isomorphic endpoints.
  The conclusion follows.
\end{proof}

Below we will use that,
by~\eqref{eq:its-locally-deRham}, 
\begin{equation}
  \label{eq:cohomology_good_freqs}
  \xi \in \Gamma_{\pmb{\omega}}
  \Longrightarrow
  \sol_{\xi \cdot \pmb{\omega}} \cong \sol_{0}
  \Longrightarrow
  \HH^{q}_\xi(\cinfty(M)) \cong H^q_{\mathrm{dR}}(M),
  \quad \forall q \in \{0, \ldots, n\}.
\end{equation}

\begin{Prop}
  \label{prop:finiteness1}
  Suppose that
  $\Gamma_{\pmb{\omega}} \neq \{0 \}$.
  Given $q \in \{0, \ldots, n\}$,
  if
  $\HH_{\dd'}^{q}(\cinfty(\Omega))$
  is finite dimensional, then
  $\HH^q_{\xi}(\cinfty(M)) = \{0\}$
  for every $\xi \in \Z^m$
  -- in particular,
  $H^q_{\mathrm{dR}}(M) = \{0\}$ --,
  hence
  $\HH_{\dd'}^{q}(\cinfty(\Omega)) = \{0\}$.
\end{Prop}
\begin{proof}
  Pick a non-zero
  $\eta \in \Gamma_{\pmb{\omega}}$
  and a non-vanishing function
  $f \in \mathscr{S}_\eta(M)$.
  Given
  $u \in \cinfty(M; \Lambda^q)$ and
  $\xi \in \Z^m$,
  we have that
  \begin{equation*}
    \dd (fu) + i ((\xi + \eta) \cdot \pmb{\omega}) \wedge (fu)
    = \dd f \wedge u + f \dd u + i f (\xi \cdot \pmb{\omega} + \eta \cdot \pmb{\omega}) \wedge u)
    = f (\dd u + i (\xi \cdot \pmb{\omega}) \wedge u)
  \end{equation*}
  hence multiplication by $f$ maps
  $\HH^q_\xi(\cinfty(M))$ to
  $\HH^q_{\xi + \eta}(\cinfty(M))$;
  since multiplication by $1/f$
  reverses this job
  we have
  $\HH^q_\xi(\cinfty(M)) \cong \HH^q_{\xi + \eta}(\cinfty(M))$
  for every
  $\xi \in \Z^m$.
  Thus,
  if some
  $\HH^q_\xi(\cinfty(M))$ is non-zero,
  then so are infinitely many of them,
  and the left-hand side
  of~\eqref{eq:into_sum}
  must contain
  infinitely many copies of it:
  injectivity of~\eqref{eq:into_sum}
  would then lead us
  to a contradiction,
  hence proving our first claim.
  Now, it follows from
  Theorem~\ref{thm:finitedimensional}
  (with $F = \eset$)
  that 
  $\HH_{\dd'}^{q}(\cinfty(\Omega))$ vanishes.
\end{proof}

\begin{Cor}
  If
  $0 < \dim \HH_{\dd'}^{q}(\cinfty(\Omega)) < \infty$
  for some
  $q$,
  then
  $\Gamma_{\pmb{\omega}} = \{0\}$.
\end{Cor}
Interesting special cases
of the results above
are obtained for
$q \in \{0, n\}$,
for then
$H^q_{\mathrm{dR}}(M)$ is one-dimensional
(compare with~\cite[Lemma~2.2]{bcp96}).
\begin{Prop}
  Suppose that:
  \begin{enumerate}
  \item $\Gamma_{\pmb{\omega}} = \Z^m$;
  \item $H^q_{\mathrm{dR}}(M) = \{0\}$; and
  \item $\dd': \cinfty(\Omega; \LLambda^{q - 1}) \to \cinfty(\Omega; \LLambda^{q})$
    has closed range.
  \end{enumerate}
  Then
  $\HH_{\dd'}^{q}(\cinfty(\Omega)) = \{0\}$.
\end{Prop}
\begin{proof}
  The first hypothesis
  ensures that all the factors
  in the direct product
  in~\eqref{eq:into_product}
  are copies of
  $H^q_{\mathrm{dR}}(M)$,
  which is zero
  by the second hypothesis;
  the third one
  implies that
  $\HH_{\dd'}^{q}(\cinfty(\Omega))$
  embeds there,
  which is therefore
  also zero.
\end{proof}
\begin{Thm}
  \label{thm:finiteness_when_mysterious}
  If~\eqref{eq:vanishing_mysterious_cohomologies}
  holds, then
  \begin{equation}
    \label{eq:finite_iff_deRham}
    \text{$\HH_{\dd'}^{q}(\cinfty(\Omega))$ is finite dimensional}
    \Longleftrightarrow
    \HH_{\dd'}^{q}(\cinfty(\Omega)) \cong H^q_{\mathrm{dR}}(M).
  \end{equation}
\end{Thm}
\begin{proof}
  We have
  by~\eqref{eq:vanishing_mysterious_cohomologies}
  and~\eqref{eq:cohomology_good_freqs}:
  \begin{equation*}
    \bigoplus_{\xi \in \Z^m} \HH_\xi^q(\cinfty(M))
    =
    \bigoplus_{\xi \in \Gamma_{\pmb{\omega}}} \HH_\xi^q(\cinfty(M))
    \cong
    \bigoplus_{\xi \in \Gamma_{\pmb{\omega}}} H_{\mathrm{dR}}^q(M)
  \end{equation*}
  which is
  isomorphic with
  $\HH_{\dd'}^{q}(\cinfty(\Omega))$
  by Theorem~\ref{thm:finitedimensional}
  provided the latter space
  is finite dimensional.
  Two possibilities arise:
  either
  $H_{\mathrm{dR}}^q(M)$
  is zero
  (in which case
  $\HH_{\dd'}^{q}(\cinfty(\Omega))$ vanishes too)
  or not;
  in the latter case,
  there must be at most
  finitely many indices
  in the last direct sum
  (by finite dimensionality)
  i.e.~$\Gamma_{\pmb{\omega}}$
  is a finite set.
  Since this is
  a subgroup of $\Z^m$,
  we have that
  $\Gamma_{\pmb{\omega}} = \{0\}$
  i.e.~there is
  a single term in that direct sum,
  namely
  $\HH^q_0(\cinfty(M)) = H_{\mathrm{dR}}^q(M)$.
  In both cases
  we get~\eqref{eq:finite_iff_deRham}.
\end{proof}

\subsection{Calculations on surfaces}
\label{sec:surfaces}

In the case
$\dim M = 2$, we know that
\begin{equation*}
  \dim H_{\mathrm{dR}}^q(M) =
  \begin{cases}
    1, &\text{if $q = 0, 2$}, \\
    2\mathsf{g}, &\text{if $q = 1$},
  \end{cases}
\end{equation*}
in which
$\mathsf{g}$ is the genus of
$M$. 
It follows from the
Atiyah-Singer Index Theorem
that the index of
the elliptic complex
$\dd'_\xi$ on $M$,
that is,
\begin{equation}
  \label{eq:indexxi}
  \sum_{q = 0}^2 (-1)^q \dim \HH_{\xi}^q(\cinfty(M)),
\end{equation}
depends only on
the principal symbol of $\dd'_\xi$,
which is the same as that of
the exterior derivative.
In particular,
\eqref{eq:indexxi} does not depend on
$\xi$, hence equals
\begin{equation*}
  \sum_{q = 0}^2 (-1)^q \dim \HH_{\xi}^q(\cinfty(M))
  = \sum_{q = 0}^2 (-1)^q \dim H_{\mathrm{dR}}^q(M)
  = \chi(M)
  = 2 - 2 \mathsf{g},
\end{equation*}
the Euler characteristic of
$M$.
We already knew this for
$\xi \in \Gamma_{\pmb{\omega}}$
by~\eqref{eq:cohomology_good_freqs};
the extra information comes for 
$\xi \notin \Gamma_{\pmb{\omega}}$,
in which case by definition
\begin{equation*}
  \HH_{\xi}^0(\cinfty(M))
  = \sol_{\xi \cdot \pmb{\omega}}(M)
  = \{0\}.
\end{equation*}
Also,
it follows from
the second part of
Lemma~\ref{lem:transposes}
that
\begin{equation*}
  \HH^2_{\xi} (\cinfty(M))
  \cong \HH^0_{-\xi} (\D'(M))^*
  \cong \HH^0_{-\xi} (\cinfty(M))^*
\end{equation*}
also vanishes, since
$-\xi \notin \Gamma_{\pmb{\omega}}$.
We conclude that for 
$\xi \notin \Gamma_{\pmb{\omega}}$
we have
\begin{equation*}
  \dim \HH_{\xi}^q(\cinfty(M)) =
  \begin{cases}
    0, &\text{if $q = 0, 2$}, \\
    2 \mathsf{g} - 2, &\text{if $q = 1$}.
  \end{cases}
\end{equation*}

\subsection*{Case $\mathsf{g} = 0$}

In the case $M$ is the $2$-sphere, we always have
$\Gamma_{\pmb{\omega}} = \Z^m$, since
every closed
$1$-form is exact by
simply connectedness.
Hence:
\begin{itemize}
\item $\HH_{\dd'}^{0}(\cinfty(\Omega))$
  and $\HH_{\dd'}^{2}(\cinfty(\Omega))$
  are infinite dimensional
  (Proposition~\ref{prop:finiteness1});
\item $\HH_{\dd'}^{1}(\cinfty(\Omega))$
  is finite dimensional
  if and only if
  $\dd': \cinfty(\Omega) \to \cinfty(\Omega; \LLambda^{1})$
  has closed range,
  in which case
  (Theorem~\ref{thm:finitedimensional})
  \begin{equation*}
    \HH_{\dd'}^{1}(\cinfty(\Omega))
    \cong
    \bigoplus_{\xi \in \Z^m} \HH_{\xi}^1(\cinfty(M))
    \cong
    \bigoplus_{\xi \in \Z^m} H_{\mathrm{dR}}^1(M)
    = \{0\}.
  \end{equation*}
\end{itemize}

\subsection*{Case $\mathsf{g} = 1$}

In the case
$M$ is the $2$-torus,
since no de Rham cohomology space vanishes,
we must have that every
$\HH_{\dd'}^{q}(\cinfty(\Omega))$
is infinite dimensional,
unless
$\Gamma_{\pmb{\omega}} = \{0\}$
(Proposition~\ref{prop:finiteness1}).
In this case we have that
\begin{equation*}
  \bigoplus_{\xi \in \Z^m} \HH_{\xi}^q(\cinfty(M))
  = H_{\mathrm{dR}}^q(M) \oplus \bigoplus_{\xi \notin \Gamma_{\pmb{\omega}}} \HH_{\xi}^q(\cinfty(M))
  = H_{\mathrm{dR}}^q(M)
\end{equation*}
is in particular finite dimensional,
hence~\eqref{eq:finite_iff_deRham} holds
for each
$q \in \{0, 1, 2\}$
(finiteness granted when
$q = 0$
thanks to
Theorem~\ref{thm:finitedimensional}).

\subsection*{Case $\mathsf{g} \geq 2$}

As in the previous case,
every
$\HH_{\dd'}^{q}(\cinfty(\Omega))$
is infinite dimensional,
except when
$\Gamma_{\pmb{\omega}} = \{0\}$,
in which case:
\begin{itemize}
\item for $q = \{0, 2\}$,
  we have
  \begin{equation*}
    \bigoplus_{\xi \in \Z^m} \HH_{\xi}^q(\cinfty(M))
    = H_{\mathrm{dR}}^q(M) \oplus \bigoplus_{\xi \notin \Gamma_{\pmb{\omega}}} \HH_{\xi}^q(\cinfty(M))
    = H_{\mathrm{dR}}^q(M)
  \end{equation*}
  hence~\eqref{eq:finite_iff_deRham} holds,
  with finiteness ensured
  at least for $q = 0$;
\item for $q = 1$,
  we have that
  \begin{equation*}
    \bigoplus_{\xi \in \Z^m} \HH_{\xi}^1(\cinfty(M))
    = H_{\mathrm{dR}}^1(M) \oplus \bigoplus_{\xi \notin \Gamma_{\pmb{\omega}}} \HH_{\xi}^1(\cinfty(M)) 
  \end{equation*}
  is infinite dimensional,
  hence so is
  $\HH_{\dd'}^{1}(\cinfty(\Omega))$
  by Proposition~\ref{prop:into_sum}.
\end{itemize}

\subsection{On the existence of an isomorphism under global solvability}
\label{exa:dm_fails}

Finally,
back again to the case when
$\omega_1, \ldots, \omega_m$ are real,
the main result in~\cite{dm16}
-- in which
$M = \TT^n$ --
essentially states that,
under a hypothesis that is equivalent to
the global solvability in degree~$1$
of the operator $\dd'$,
we have
\begin{equation}
  \label{eq:dm_isoparaotoro1}
  \HH^{q}_{\dd'}(\cinfty(\TT^{m+n}))
  \cong
  \cinfty(\TT^r) \otimes H^q_{\mathrm{dR}}(\TT^m).
\end{equation}
This inspired us to prove our
isomorphism~\eqref{IsomorphismSobreGammaomega},
which considers only the cohomology of
$q$-forms associated with the cluster
$\Gamma_{\pmb{\omega}}$.
We will show below
that an isomorphism similar
to~\eqref{eq:dm_isoparaotoro1}
is not true for general
$M$,
even assuming solvability
in the first degree.

Indeed, back to the case
$\dim M = 2$,
$\mathsf{g} \geq 2$,
let
$\pmb{\omega} \dfn \{ \lambda \vartheta_1 \}$
(corank~$1$)
in which
$\vartheta_1$ is as in the proof of
Proposition~\ref{Prop:non-simult-appr-charac}
($d = 2 \mathsf{g} > 0$ there).
Therefore, if
$\lambda \in \R \setminus \Q$,
then
$\Gamma_{\pmb{\omega}} = \{0\}$
thus $r = 0$,
hence
$\cinfty(\TT^r) \otimes H^1_{\mathrm{dR}}(M) \cong H^1_{\mathrm{dR}}(M)$
cannot be isomorphic with the
infinite dimensional
$\HH_{\dd'}^{1}(\cinfty(\Omega))$,
even when
$\dd': \cinfty(\Omega) \to \cinfty(\Omega; \LLambda^1)$
has closed range
(for instance, if
$\lambda$ is a non-Liouville number,
by Theorem~\ref{Thm:solv-diop}).

\def\cprime{$'$}

\end{document}